\documentclass[final]{siamltex}

\usepackage{amssymb,amsmath}
\usepackage{algorithm,algorithmic}
\usepackage{epsfig,epstopdf}

\def\etal{\emph{et al.}}
\newcommand{\remove}[1]{}
\newcommand{\chola}{Cholesky decomposition algorithm}
\newcommand{\cholas}{Cholesky decomposition algorithms}
\newcommand{\chold}{Cholesky decomposition}
\newcommand{\cholf}{Cholesky factorization}
\newcommand{\scalapack}{$ScaLAPACK$}
\newcommand{\lapack}{$LAPACK$}
\newcommand{\lt}{\left}
\newcommand{\rt}{\right}

\newcommand{\ignore}[1]{}

\renewcommand{\algorithmicforall}{\textbf{parallel for each}}
\newcommand{\yes}{$\checkmark$}
\newcommand{\no}{$\times$}
\renewcommand{\algorithmicrequire}{\textbf{Input:}}
\renewcommand{\algorithmicensure}{\textbf{Output:}}

\title{Communication-Optimal Parallel and Sequential Cholesky Decomposition
\thanks{A preliminary version of this paper was accepted to
the 21st ACM Symposium on Parallelism in Algorithms and Architectures (SPAA'09) \cite{BDHS09-SPAA}.}}
\author{
Grey Ballard
\thanks{Computer Science Department,
University of California, Berkeley, CA 94720.
Research supported by Microsoft and Intel funding (Award $\#$20080469) and by matching funding by U.C. Discovery (Award $\#$DIG07-10227).
(ballard@eecs.berkeley.edu).} \and
James Demmel
\thanks{Mathematics Department and CS Division,
University of California, Berkeley, CA 94720.
Research supported by Microsoft and Intel funding (Award $\#$20080469) and by matching funding by U.C. Discovery (Award $\#$DIG07-10227).
(demmel@cs.berkeley.edu).} \and
Olga Holtz
\thanks{Departments of Mathematics,
University of California, Berkeley and
Technische Universit\"at Berlin. O. Holtz acknowledges support of
the Sofja Kovalevskaja programme of Alexander von Humboldt
Foundation.
(holtz@math.berkeley.edu).} \and
Oded Schwartz
\thanks{Departments of Mathematics,
Technische Universit\"at Berlin, 10623 Berlin, Germany.
This work was done while visiting University of California, Berkeley.
(odedsc@math.tu-berlin.de). }}

\begin{document}
\maketitle
\begin{abstract}
Numerical algorithms have two kinds of costs: arithmetic and communication, by which we mean either moving data between levels of a memory hierarchy (in the sequential case) or over a network connecting processors (in the parallel case). Communication costs often dominate arithmetic costs, so it is of interest to design algorithms minimizing communication. In this paper we first extend known lower bounds on the communication cost (both for bandwidth and for latency) of conventional ($O(n^3)$) matrix multiplication to \cholf, which is used for solving dense symmetric positive definite linear systems. Second, we compare the costs of various \chold\ implementations to these lower bounds and identify the algorithms and data structures that attain them.  In the sequential case, we consider both the two-level and hierarchical memory models.  Combined with prior results in \cite{DemmelGrigoriHoemmenLangou08a,DemmelGrigoriHoemmenLangou08b,DemmelGrigoriXiang08}, this gives a set of communication-optimal algorithms for $O(n^3)$ implementations of the three basic factorizations of dense linear algebra: LU with pivoting, QR and Cholesky.  But it goes beyond this prior work on sequential LU by optimizing communication for any number of levels of memory hierarchy.
\end{abstract}

\thispagestyle{empty}

\setcounter{page}{1} \pagestyle{plain} \pagenumbering{arabic}

%% update according to publishing house %%  \category{F.2.1}{Theory of Computation}{Analysis of Algorithms and Problem Complexity}[Numerical Algorithms and Problems,Computations on matrices]
%% update according to publishing house %%  \terms{Algorithms, Performance.}


\pagestyle{myheadings}
\thispagestyle{plain}
\markboth{BALLARD, DEMMEL, HOLTZ, AND SCHWARTZ}{COMMUNICATION-OPTIMAL CHOLESKY DECOMPOSITION}

\section{Introduction}

Let $A$ be a real symmetric and positive definite matrix.
Then there exists a real lower triangular matrix $L$ so
that $A=L \cdot L^T$ ($L$ is unique if we restrict its diagonal
elements to be positive).
This is called the \chold.
We are interested in finding efficient parallel and sequential algorithms
for the \chold. Efficiency is measured both by
the number of arithmetic operations and by the amount of
communication, either between levels of a memory hierarchy on
a sequential machine, or between processors communicating over a
network on a parallel machine.
Since the time to move one word of data typically exceeds the
time to perform one arithmetic operation by a large and growing
factor \cite{GrahamSnirPatterson04}, our goal will be to minimize communication.

\subsection{Communication model}

We model communication costs in more detail as follows. In the
sequential case, with two levels of memory hierarchy (fast and
slow), communication means reading data items ({\em words}) from
slow memory to fast memory and writing data from fast memory to
slow memory. Words that are contiguous can be read or
written in a bundle which we will call a {\em message}. We assume
that a message of $n$ words can be communicated between fast and
slow memory in time $\alpha+\beta n$ where $\alpha$ is the {\em
latency} (seconds per message) and $\beta$ is the {\em inverse
bandwidth} (seconds per word). We define the {\em bandwidth cost} of an algorithm to be the total number of words communicated and the {\em latency cost} of an algorithm to be the total number of messages communicated.  We assume that the matrix being factored initially resides in slow memory, and is too large to fit
in the smaller fast memory. Our goal then is to minimize the total
number of words and the total number of messages communicated
between fast and slow memory.\footnote{The sequential
communication model used here is sometimes called the
\emph{two-level I/O model} or \emph{disk access machine (DAM)}
model (see \cite{AggarwalVitter88,BenderBrodalFagerbergJacobVicari07,ChowdhuryRamachandran06}).  Our bandwidth cost model follows that of
\cite{HongKung81} and \cite{IronyToledoTiskin04} in that it
assumes the block-transfer size is one word of data ($B=1$ in the
common notation).  However, our model allows message sizes to vary from one word up to the maximum number of words that can fit in fast memory.}

In the parallel case, we are interested in the communication among
the processors.  As in the sequential case, we assume that a message
of $n$ consecutively stored words can be communicated in time
$\alpha+\beta n$. This cost includes the time required to ``pack'' non-contiguous words into a single message, if necessary.  We assume that the matrix is initially evenly
distributed among all $P$ processors, and that there is only enough
memory to store about $1/P$-th of a matrix per processor.
As before, our goal is to minimize the number of words and messages
communicated.  In order to measure  the communication complexity of
a parallel algorithm, we will count the number of words and messages
communicated along the critical path of the algorithm.  That is, if many pairs of processors send identical messages simultaneously, the cost along the critical path of the algorithm is that of one message (we do not consider communication resource contention among processors).

\subsection{Communication Lower Bounds}

We consider {\em classical} algorithms for \chold,
i.e., those that perform ``the usual'' $O(n^3)$ arithmetic operations,
possibly reordered by associativity and commutativity of addition.
That is, our results do not apply when using
distributivity to reorganize the algorithm
(such as Strassen-like algorithms);
we also assume no pivoting is performed.
We define ``classical'' more carefully later.
We show that the communication complexity of any such Cholesky algorithm
shares essentially the same lower bound as does the
classical matrix multiplication.

\begin{theorem}[Main Theorem] \label{thm:main}
Any sequential or parallel classical algorithm for the \chold\ of
$n$-by-$n$ matrices can be transformed into a classical algorithm
for $\frac{n}{3}$-by-$\frac{n}{3}$ matrix multiplication, in such
a way that the bandwidth cost of the matrix multiplication algorithm is
at most a constant times the bandwidth cost of the Cholesky algorithm.
\end{theorem}

Therefore any bandwidth cost lower bound for classical matrix multiplication applies to classical Cholesky decomposition, in a Big-O sense.  In particular, since a sequential classical $n$-by-$n$ matrix multiplication algorithm has a bandwidth cost lower bound of $\Omega ( n^3 / M^{1/2} )$ where $M$ is the fast memory size \cite{HongKung81,IronyToledoTiskin04}, classical Cholesky decomposition has the same lower bound (we discuss the parallel case later).

To get a latency cost lower bound, we use the simple observation
\cite{DemmelGrigoriHoemmenLangou08a} that the number
of messages is at least the bandwidth cost lower bound divided by
the maximum message size, and that the maximum message size
is at most fast memory size in the sequential case (or the
local memory size in the parallel case). So for sequential
matrix multiplication this means the latency cost lower bound is
$\Omega ( n^3 / M^{3/2} )$.

\subsection{Communication Upper Bounds}

\subsubsection{Two-level memory model}

In the case of matrix multiplication on a sequential machine, a well-known algorithm attains the bandwidth cost lower bound \cite{HongKung81}.  It computes $C = A \cdot B$ by performing blocked matrix multiplication that multiplies and accumulates $\sqrt{\frac{M}{3}}$-by-$\sqrt{\frac{M}{3}}$ submatrices of $A$, $B$ and $C$. If each matrix is stored so that the $\frac{M}{3}$
entries of each of its submatrices are contiguous
(not the case with columnwise or rowwise storage)
the latency cost lower bound is reached; we call such
a data structure {\em block-contiguous storage} and describe
it in more detail below.
Alternatively, one could try to copy $A$ and $B$ from their input
format (say columnwise) to contiguous block storage doing
(asymptotically) no more communication than the subsequent
matrix multiplication; we will see this is possible provided
$M = \Omega(n)$.
There will be analogous requirements for Cholesky to attain
its latency cost lower bound.

In particular, we draw the following conclusions about the communication costs of sequential
classical Cholesky decomposition in the two-level memory model, as summarized in Table~\ref{tbl:sequential}:
\newcounter{saveenum}
\begin{enumerate}
    \item \label{con:naive} ``Na{\"\i}ve'' sequential variants of Cholesky that operate on single rows and columns (be they left-looking, right-looking, etc.) attain neither the bandwidth nor the latency cost lower bound.

    \item \label{con:known-bandwidth} A sequential blocked algorithm used in LAPACK (with the correct block size) attains the bandwidth cost lower bound, as do the recursive algorithms in \cite{AhmedPingali00,AndersenGustavsonWasniewski01,Gustavson97,GustavsonJonsson01,SimecekTvrdik04}. A recursive algorithm analogous to Toledo's LU algorithm \cite{Toledo97} attains the bandwidth cost lower bound in nearly all cases, expect possibly for an $O(\log n)$ factor in the narrow range $\frac{n^2}{\log ^2 n} < M < n^2$.

    \item \label{con:known-latency}  Whether the LAPACK algorithm also attains the latency cost lower bound depends on the matrix layout: If the input matrix is given in row-wise or column-wise format, and this is not changed by the algorithm, then the latency cost lower bound cannot be attained. But if the input matrix is given in contiguous block storage, or $M = \Omega(n)$ so that it can be copied quickly to contiguous block format, then the latency cost lower bound can be attained by the \lapack\ algorithm. Toledo's algorithm cannot minimize latency (at least when $M > n^{2/3}$)  See also Conclusion~\ref{con:new-recursive} below.

\setcounter{saveenum}{\value{enumi}}
\end{enumerate}

\subsubsection{Hierarchical memory model}

Since most computers have multiple levels (e.g., L1, L2, L3 caches, main memory, and disk), we also consider the communication costs in the hierarchical memory model.
In this case, an optimal algorithm should simultaneously
minimize communication between {\em all} pairs of
adjacent levels of memory hierarchy (e.g., minimize
bandwidth and latency costs between L1 and L2, between L2 and L3, etc.).

In the case of sequential matrix multiplication, bandwidth cost is minimized
in this sense by simply applying the usual blocked algorithm
recursively, where each level of recursion multiplies
matrices that fit at a particular level of the memory hierarchy,
by using the blocked algorithm to multiply submatrices that
fit in the next smaller level. This is easy since matrix
multiplication naturally breaks into smaller matrix multiplications.

For matrix multiplication to minimize latency across all memory
hierarchy levels, it is necessary for all submatrices of all sizes
to be stored contiguously. This leads to a data structure
variously referred to as {\em recursive block storage}, {\em Morton ordering}, or storage
using {\em space-filling curves}, and described in
\cite{AhmedPingali00,Bader08,ElmrothGustavsonJonssonKagstrom04,FrigoLeisersonProkopRamachandran99,Wise00}.

Finally, sequential matrix multiplication can achieve
communication optimality as just described in one of two ways: (1)
We can choose the number of recursion levels and sizes of the
subblocks with prior knowledge of the number of levels and sizes
of the levels of memory hierarchy, a {\em cache-aware}
process called tuning. (2) We carry out the recursion down to
1-by-1 blocks (or some other small constant size), repeatedly
dividing block sizes by 2 (perhaps padding submatrices to have
even dimensions as needed). Such an algorithm is called {\em
cache-oblivious} \cite{FrigoLeisersonProkopRamachandran99}, and
has the advantage of simplicity and portability compared to a
cache-aware algorithm, though it might also have more
overhead in practice.

It is indeed possible for sequential Cholesky to be organized to
be optimal across multiple memory hierarchy levels in all the senses
just described, assuming we use recursive block storage:
\begin{enumerate}
\setcounter{enumi}{\value{saveenum}}
    \item \label{con:known-cache-oblivious} The recursive algorithm modeled on Toledo's LU can be implemented in a cache-oblivious way so as to minimize bandwidth cost, but not latency cost.\protect\footnote{Toledo's algorithm is designed to retain numerical stability for the LU case. The algorithm of \cite{Gustavson97} deals with the Cholesky case only and therefore requires no pivoting for numerical stability. Thus a simpler recursion suffices, and the latency cost diminishes.}

    \item \label{con:new-recursive} The cache-oblivious recursive Cholesky algorithm which first appeared in \cite{Gustavson97} (and can also be found in \cite{AhmedPingali00,AndersenGustavsonWasniewski01,GustavsonJonsson01,SimecekTvrdik04}) minimizes both bandwidth and latency costs (assuming recursive block storage) for all matrices across all memory hierarchy levels. None of the other algorithms have this property.

\setcounter{saveenum}{\value{enumi}}
\end{enumerate}

\subsubsection{Parallel model}

Finally, we address the case of parallel Cholesky,
where there are $P$ processors connected by a network
with latency $\alpha$ and reciprocal bandwidth $\beta$.
We consider here only the memory-scalable case, where
each processor's local memory is of size $M = O(n^2/P)$,
so that only $O(1)$ copies of the matrix are stored overall
(the so-called ``2D case'').  See \cite{IronyToledoTiskin04}
for the general lower bound, including 3D, for matrix multiplication.  A 3D algorithm for LU decomposition without pivoting (along with an algorithm for triangular solve) is presented with communication analysis in \cite{IronyToledo01}.  While the algorithms minimize the total volume of communication, they do not minimize communication costs along the critical path as measured in this paper.

For the 2D case, the consequence of our Main Theorem is again a bandwidth cost
lower bound of the form
$\Omega (n^3 / (P M^{1/2})) = \Omega (n^2 / P^{1/2})$,
and a latency cost lower bound of the form
$\Omega (n^3 / (P M^{3/2})) = \Omega (P^{1/2})$.

\begin{enumerate}
\setcounter{enumi}{\value{saveenum}}
    \item \label{con:parallel} The Cholesky algorithm in \scalapack\ \cite{SCALAPACK} attains a matching upper bound. It does so by partitioning the matrix into submatrices and distributing them to the processors in a block cyclic manner. With the right choice of block size $b$, namely $b = \Theta\lt(\frac{n}{\sqrt{P}}\rt)$, it attains the above bandwidth and latency cost lower bounds to within a factor of $\log P$. This is summarized in Table~\ref{tbl:parallel}.  See \S~\ref{sec:scalapack} for details.
\end{enumerate}

\begin{table}[!ht]
\footnotesize
  \centering
\begin{tabular}{|ll||c|c|c|} \hline
& & & & Cache \\
& & Bandwidth & Latency  & Oblivious\\
\hline \hline
\multicolumn{2}{|l||}{Lower bound}
& $\Omega\lt(\frac{n^3}{\sqrt{M}}\rt)$
&  $\Omega\lt(\frac{n^3}{M^{3/2}}\rt)$ &\\
\hline \hline
\multicolumn{2}{|l||}{Na{\"\i}ve: left/right looking} & & & \\
& Column-major
&  $\Theta(n^3)$ &  $\Theta\lt(  n^2+\frac{n^3}{M}   \rt)$
& \yes\\
\hline
\multicolumn{2}{|l||}{\lapack\ \cite{LAPACK}} & & & \\
& Column-major
&  $O\lt(\frac{n^3}{\sqrt{M}} \rt)$
& $O\lt(\frac{n^3}{{M}} \rt)$
& \no\\
& Block-contiguous
&  $O\lt(\frac{n^3}{\sqrt{M}}\rt)$
& $O\lt(\frac{n^3}{{M^{3/2}}}\rt)$
& \no \\
\hline
\multicolumn{2}{|l||}{Rectangular Recursive \cite{Toledo97}}  & & & \\
& Column-major
& $\Theta\lt(\frac{n^3}{\sqrt{M}} + n^2 \log n \rt)$
& $\Omega \left( \frac{n^3}{M} \right)$
& \yes \\
& Block-contiguous
& $\Theta\lt(\frac{n^3}{\sqrt{M}} + n^2 \log n \rt)$
& $\Omega \left( n^2 \right)$ & \yes\\
\hline
\multicolumn{2}{|l||}{Square Recursive \cite{AhmedPingali00,AndersenGustavsonWasniewski01,Gustavson97,GustavsonJonsson01,SimecekTvrdik04} }  & & & \\
& ``Recursive Packed Format''   \cite{AndersenGustavsonWasniewski01}
& $O\lt(\frac{n^3}{\sqrt{M}} \rt)$
& $\Omega \left( \frac{n^3}{M} \right)$
& \yes \\
& Column-major  \cite{Gustavson97}
& $O\lt(\frac{n^3}{\sqrt{M}}\rt)$
& $O\lt(\frac{n^3}{M}\rt)$
& \yes\\
& Block-contiguous \cite{AhmedPingali00}
& $O\lt(\frac{n^3}{\sqrt{M}}\rt)$
& $O \left( \frac{ n^3}{M^{3/2}} \right)$
& \yes \\
\hline
\end{tabular}
  \protect\caption{Sequential bandwidth and latency costs: lower bound vs. algorithms.
           $M$ denotes the size of the fast memory. We assume $n^2>M$ in this table. FLOPs count of all is
           $O(n^3)$. Refer to \S \ref{sec:upperbound} for
           definitions of terms and details.
  }
\label{tbl:sequential}
\end{table}
\bigskip

\begin{table}[!ht]
\footnotesize
  \centering
\begin{tabular}{|l||c|c|c|} \hline
& Bandwidth & Latency & FLOPS \\
\hline \hline
Lower-bound & & &\\
~~~General & $\Omega\lt(\frac{n^3}{P\sqrt{M}}\rt)$ &  $\Omega\lt(\frac{n^3}{PM^{3/2}}\rt)$ & $\Omega\lt(\frac{n^3}{P}\rt)$ \\
~~~2D layout: \; $M=O\lt(\frac{n^2}{P}\rt)$ & $\Omega\lt(\frac{n^2}{\sqrt{P}}\rt)$ &  $\Omega\lt( \sqrt{P} \rt)$ & $\Omega\lt(\frac{n^3}{P}\rt)$ \\
\hline \hline
\scalapack\ \cite{SCALAPACK}  & & &\\
~~~General & $O\lt(\lt(\frac{n^2}{\sqrt{P}}+nb\rt)\log P\rt)$  & $O\lt(\frac{n}{b}\log P\rt) $   & \\
~~~Choosing $b=\Theta\lt(\frac{n}{\sqrt{P}}\rt)$ & $O\lt(\frac{n^2}{\sqrt{P}}\log P\rt)$ & $O(\sqrt{P}\log P) $ & $O\lt(\frac{n^3}{P}\rt)$ \\
\hline
\end{tabular}

  \protect\caption{Parallel bandwidth and latency costs: lower bound vs. algorithms.
    $M$ denotes the size of the memory of each processor. $P$ is the number of processors. $b$ is the block size.
    Refer to \S \ref{sec:upperbound} for
           definitions of terms and details.
  }\label{tbl:parallel}
\end{table}

The rest of this paper is organized as follows. In \S \ref{sec:lowerbound} we show the reduction from
matrix multiplication to \chold, thus extending the bandwidth cost
lower bounds of \cite{HongKung81} and \cite{IronyToledoTiskin04}
to a bandwidth cost lower bound for the sequential and parallel
implementations of \chold. We also discuss latency cost lower bounds.
In \S \ref{sec:upperbound} we present known \cholas\ and
compare their bandwidth and latency costs with the lower bounds.

A shorter version of this paper (an extended abstract) appeared in the proceedings of the SPAA 2009 conference; the shorter version includes the reduction proof of the lower bounds and Tables~\ref{tbl:sequential} and \ref{tbl:parallel}.  This paper proves the claims made in the tables by presenting the algorithms and their analyses in \S \ref{sec:upperbound}.

\section{Communication Lower Bounds}
\label{sec:lowerbound}

Consider an algorithm for a parallel computer  with $P$ processors
that multiplies matrices  in the ``classical'' way (the usual $O(n^3)$
arithmetic operations possibly reordered using associativity
and commutativity of addition) and each of the processors has memory of size
$M$. Irony \etal\/  \cite{IronyToledoTiskin04} showed
that at least one of the processors has to send or receive this minimal
 number of words:
\begin{theorem}[{\cite{IronyToledoTiskin04}}]
\label{thm:IronyToledoTiskin04}
Any ``classical'' implementation of matrix multiplication of $n{\times}n$ matrices $A$ and $B$ on a $P$
processor machine, each equipped with memory of size $M$, requires
that one of the processors sends or receives
 at least $\frac{n^3}{2\sqrt 2 PM^{\frac12}}-M$ words.

If $A$ and $B$ are of size $n{\times}m$ and $m{\times}r$ respectively,
then the corresponding bound is
$\frac{nmr}{2\sqrt 2 PM^{\frac12}}-M$.
\end{theorem}

As any processor has memory of size $M$, any message it sends or
receives may deliver at most $M$ words. Therefore we deduce the
following:

\begin{corollary}
\label{cor:IronyToledoTiskin04}
Any ``classical'' implementation of matrix multiplication of $n{\times}n$ matrices $A$ and $B$ on a $P$
processor machine, each equipped with memory of size $M$, requires
that one of the processors sends or receives
 at least $\frac{n^3}{2\sqrt 2 PM^{\frac32}}-1$ messages.

 If $A$ and $B$ are of size $n{\times}m$ and $m{\times}r$ respectively,
then the corresponding bound is
$\frac{nmr}{2\sqrt 2 PM^{\frac32}}-1$.
\end{corollary}

For the case of $P=1$ these coincide with the lower bounds for the bandwidth and the latency costs of the sequential case.
The lower bound on bandwidth cost of sequential matrix multiplication was previously shown (up to some multiplicative constant factor) by Hong and Kung \cite{HongKung81}.

One means of extending these lower bounds to more algorithms is by reduction.  It is easy to reduce matrix multiplication to LU decomposition of a
slightly larger order, as
the following identity shows:
\begin{equation*}
\begin{pmatrix}
I & 0 & -B \\
A & I & 0 \\
0 & 0 & I \\
\end{pmatrix}
=
\begin{pmatrix}
I &   &   \\
A & I &   \\
0 & 0 & I \\
\end{pmatrix}
\cdot
\begin{pmatrix}
I & 0 & -B \\
  & I & A \cdot B \\
  &   & I \\
\end{pmatrix}
\end{equation*}
This identity means that LU factorization  without pivoting can be used to perform
matrix multiplication; to accommodate pivoting $A$ and/or $B$ can
be scaled down to be too small to be chosen as pivots, and $A
\cdot B$ can be scaled up accordingly. Thus an $O(n^3)$
implementation of LU decomposition that only uses associativity and
commutativity of addition to reorganize its operations (thus
eliminating Strassen-like algorithms) must perform at least as
much communication as a correspondingly reorganized implementation
of $O(n^3)$ matrix multiplication.

We wish to mimic this lower bound construction for Cholesky.
Consider the following reduction from matrix multiplication to
\chold. Let $T$ be the matrix defined below, composed of 9 square
blocks each of dimension $n$; then the \chold\ of $T$ is:
\small
\begin{equation*}
T \equiv
\begin{pmatrix}
  I & A^T & -B \\
  A & I+A \cdot A^T & 0 \\
  -B^T & 0 & D
\end{pmatrix}
=
\begin{pmatrix}
  I &   &   \\
  A & I &   \\
  -B^T & (A \cdot B)^T & X
\end{pmatrix}
\cdot
\begin{pmatrix}
  I & A^T & -B \\
    & I & A \cdot B \\
    &   & X^T
\end{pmatrix}
\equiv L \cdot L^T
\end{equation*}
\normalsize
where $X$ is the Cholesky factor of $D' \equiv D - B^TB - B^TA^TAB$,
and $D$ can be any symmetric matrix such that $D'$ is positive
definite.

Thus $A \cdot B$ is computed via this \chold. Intuitively this
seems to show that the communication complexity needed for
computing matrix multiplication is a lower bound to that of
computing the \chold\ (of matrices three times larger) as $A \cdot B$
appears in $L^T$, the decomposition of $T$. Note however that $A
\cdot A^T$ also appears in $T$.  Conceivably, computing $A \cdot B$ from $A$ and $B$
incurs less communication cost if we are also given $A \cdot
A^T$, so the above is not a sufficient reduction from matrix multiplication to \chold.\protect\footnote{Note that computing $A \cdot A^T$ is
asymptotically as hard as matrix multiplication: take $A=[X , 0 ;
Y^T , 0]$. Then $A \cdot A^T = [ *, XY ;
* ,* ]$} Let us instead consider the following approach to prove the lower bound.

In addition to the real numbers $\mathbb{R}$,
consider new ``starred'' numerical quantities, called $1^*$ and
$0^*$, with arithmetic properties detailed in the following tables.
$1^*$ and $0^*$ mask any real value in addition/substraction
operation, but behave similarly to $1 \in \mathbb{R}$ and
$0 \in \mathbb{R}$ in multiplication and division operations.

\begin{table}[ht]
  \centering

\begin{tabular}{c c c} % Ordering of the tables begins

\begin{tabular}{|c||c|c|c|} % Addition Table begins
\hline
  $\pm$ & $1^*$ & $0^*$ & $y $ \\
\hline \hline
  $1^*$ & $1^*$ & $1^*$ & $1^*$ \\
  \hline
  $0^*$ & $1^*$ & $0^*$ & $0^*$ \\
  \hline
  $x$ & $1^*$ & $0^*$ & $x \pm y$ \\
\hline
\end{tabular} % Addition Table ends

&
\begin{tabular}{|c||c|c|c|} % Multiplication Table begins
\hline
  $\cdot$ & $1^*$ & $0^*$ & $y$ \\
\hline \hline
  $1^*$ & $1^*$ & $0^*$ & $y$ \\
  \hline
  $0^*$ & $0^*$ & $0$ & $0$ \\
  \hline
  $x$ & $x$ & $0$ & $x \cdot y$ \\
\hline
\end{tabular} % Multiplication Table ends

 &

\begin{tabular}{|c||c|c|c|} % Division Table begins
\hline
  $/$ & $1^*$ & $0^*$ & $y \neq 0$ \\
\hline \hline
  $1^*$ & $1^*$ & $-$ & $1/y$ \\
  \hline
  $0^*$ & $0^*$ & $-$ & $0$ \\
  \hline
  $x$ & $x$ & $-$ & $x / y$ \\
\hline
\end{tabular} % Division Table ends

\\ & & \\ &

\begin{tabular}{|c||c|} % Square-root Table begins
\hline
  $\sqrt{.}$ & \\
\hline \hline
  $1^*$ & $1^*$ \\
  \hline
  $0^*$ & $0^*$ \\
  \hline
  $x \geq 0$ & $\sqrt{x}$ \\
\hline
\end{tabular} % Square-root Table ends

& \\
\end{tabular} % Ordering of the tables ends

  \protect\caption{Arithmetic Operations: $x,y$ stand for any {\em real} values. For consistency, $-0^* \equiv 0^*$ and $-1^* \equiv  1^*$.}
  \label{tbl:arith}
\end{table}

Consider this set of values and arithmetic operations.
\begin{itemize}
\item
The set is commutative with respect to addition and to multiplication (by
the symmetries of the corresponding tables).
\item
The set is associative with respect to addition: regardless of ordering
of summation, the sum is $1^*$ if one of the summands is $1^*$,
otherwise it is $0^*$ if one of the summands is $0^*$.
\item
The set is also associative with respect to multiplication: $(a
\cdot b) \cdot c = a \cdot (b \cdot c )$. This is trivial if all
factors are in $\mathbb{R}$. As $1^*$ is a multiplicative identity,
it is also immediate if some of the
factors equal $1^*$. Otherwise, at least one of the factors is $0^*$,
and the product is 0.
\item
Distributivity, however, does not hold: $1 \cdot (1^* + 1^*)= 1
\neq 2 = (1 \cdot 1^*) + (1 \cdot 1^*)$
\end{itemize}

Let us return to the construction. We set $T'$ to be:
$$
T'\equiv
\begin{pmatrix}
  I & A^T & -B \\
  A & C & 0 \\
  -B^T & 0 & C
\end{pmatrix}
$$
where $C$ has $1^*$ on the diagonal and $0^*$ everywhere else:
$$ C\equiv
\begin{pmatrix}
  1^* & 0^* & & \cdots & 0^*\\
  0^*\ & 1^* & 0^* & & \vdots \\
    & &    & \ddots & \\
  \vdots & &  &   & 0^*\\
    0^* & \cdots &  & 0^* & 1^*
\end{pmatrix}
$$

One can verify that the (unique) \chold\ of $C$
is\protect\footnote{By writing $X \cdot Y$ we mean the resulting
matrix assuming the straightforward $n^3$ matrix multiplication
algorithm. This has to be stated clearly, as the distributivity
does not hold for the starred values.}

\begin{equation}\label{eq:C}
C=
\begin{pmatrix}
  1^*    & 0      & \ldots & 0      \\
  0^*    & 1^*    &        & \vdots \\
  \vdots &        & \ddots & 0      \\
  0^*    & \cdots & 0^*    & 1^*
\end{pmatrix}
\cdot
\begin{pmatrix}
  1^*    & 0^*    & \cdots & 0^*    \\
         & \ddots & \ddots & \vdots \\
  \vdots &        &  1^*   & 0^*    \\
  0      & \hdots &        & 1^*
\end{pmatrix}
\equiv C' \cdot C'^T
\end{equation}

Note that if a matrix $X$ does not contain any ``starred'' values $0^*$ and $1^*$ then $X=C \cdot X=X \cdot C=C' \cdot X=X \cdot C'=C'^T \cdot X = X \cdot C'^T$ and $C+X=C$. Therefore, one can confirm that the \chold\ of $T'$ is:
\small
\begin{equation}
\label{eq:GEMM-to-Chol-starred}
T'\equiv
\begin{pmatrix}
  I & A^T & -B \\
  A & C & 0 \\
  -B^T & 0 & C
\end{pmatrix}
=
\begin{pmatrix}
  I &   &   \\
  A & C' &   \\
  -B^T & (A \cdot B)^T & C'
\end{pmatrix}\cdot
\begin{pmatrix}
  I & A^T & -B \\
    & C'^T & A \cdot B \\
    &   & C'^T
\end{pmatrix}
\equiv L \cdot L^T
\end{equation}
\normalsize

One can think of $C$ as masking the $A \cdot A^T$ previously
appearing in the central block of $T$, therefore allowing the lower
bound of computing $A \cdot B$ to be accounted for by the \chold, and
not by the computation of $A \cdot A^T$. More formally, let $Alg$
be any classical algorithm for \cholf. We convert it to a
matrix multiplication algorithm as follows:

\begin{algorithm}
\protect\caption{Matrix Multiplication by Cholesky Decomposition}
\label{alg:Chol-to-MM}
\begin{algorithmic}[1]
\REQUIRE Two $n{\times}n$ matrices, $A$ and $B$. \STATE Let $Alg'$
be $Alg$ updated to correctly handle the new $0^*, 1^*$ values.
\COMMENT{note that $Alg'$ can be constructed off-line.}
\STATE Construct $T'$ as in Equation (\ref{eq:GEMM-to-Chol-starred})
\STATE $L=Alg'(T')$
\RETURN $(L_{32})^T$
\end{algorithmic}
\end{algorithm}

The simplest conceptual way to do step (1) is to attach an extra
bit to every numerical value, indicating whether it is ``starred''
or not, and modify every arithmetic operation to check
this bit before performing an operation. This
increases the bandwidth cost by at most a constant factor.
Alternatively,
we can use Signalling NaNs as defined in the IEEE Floating Point
Standard \cite{IEEE08} to encode $1^*$ and $0^*$ with no extra bits.

If the instructions
implementing Cholesky are scheduled deterministically, there is
another alternative.  One can run the algorithm ``symbolically'',
propagating $0^*$ and $1^*$ arguments from the inputs forward,
simplifying or eliminating arithmetic operations whose inputs
contain $0^*$ or $1^*$. One can also eliminate operations for which
there is no path in the directed acyclic graph
(describing how outputs of each operation propagate to inputs of other
operations) to the desired output $A \cdot B$. The resulting $Alg'$
performs a strict subset of the arithmetic and memory operations of the
original Cholesky algorithm.

We note that updating $Alg$ to form $Alg'$ is done off-line,
so that step (1) does not actually take any time to perform
when Algorithm~\ref{alg:Chol-to-MM} is called.

\paragraph{Classical Cholesky Decomposition Algorithms}
We next verify the correctness of this reduction: that the output
of this procedure on input $A,B$ is indeed the multiplication $A \cdot
B$, as long as $Alg$ is a classical algorithm, in a sense
we now define carefully.

\begin{figure}[ht]
\begin{center}
\scalebox{.5}{\includegraphics{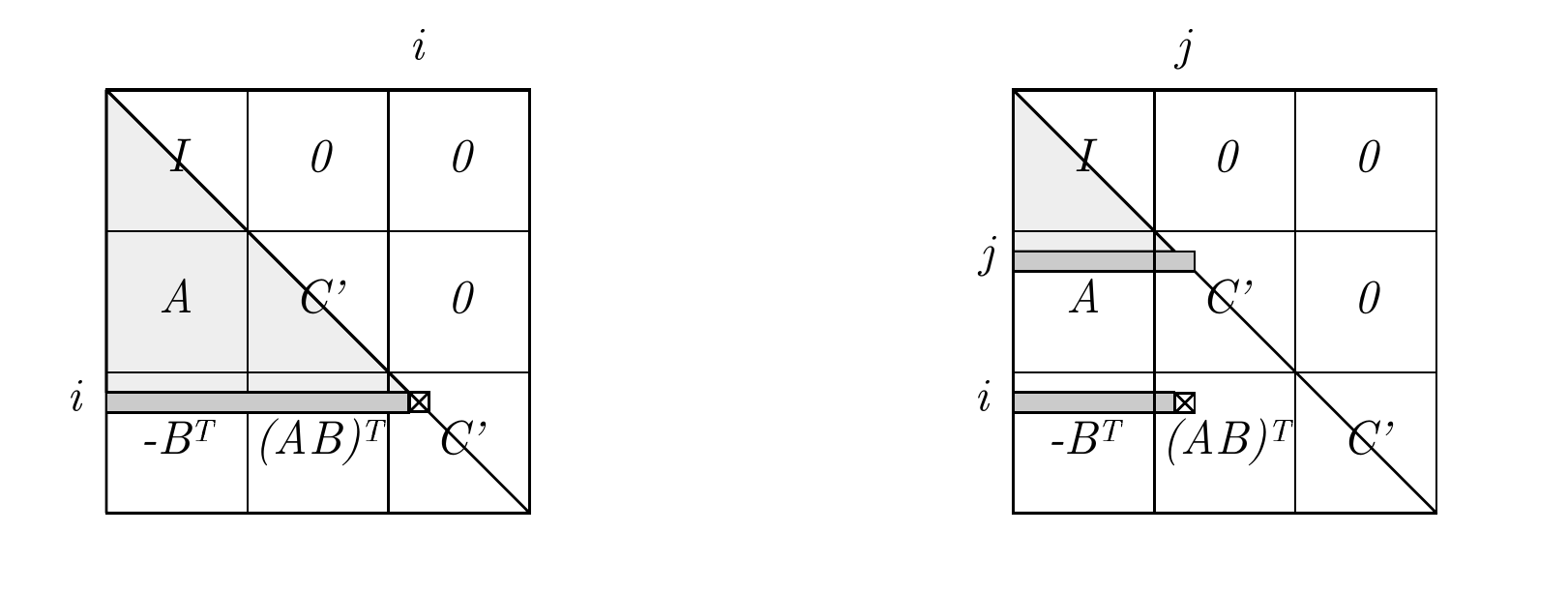}}
\protect\caption{Dependencies of $L(i,j)$, for diagonal entries
(left) and other entries (right).
\newline Dark grey represents the sets $S_{i,i}$ (left) and $S_{i,j}$ (right).
Light grey represents indirect dependencies.
 }\label{fig:dependencies}
\end{center}
\end{figure}

Let $T'=L  \cdot  L^T$ be the \chold\ of $T'$. Then we have
the following formulas:

\begin{equation}\label{eq:chol-dep-diag}
L(i,i)=\sqrt{T'(i,i)-\sum_{k \in [i-1]}(L(i,k))^2}
\end{equation}
\begin{equation}\label{eq:chol-dep}
L(i,j)=\frac{1}{L(j,j)}\left(T'(i,j)-\sum_{k\in [j-1]}L(i,k)
\cdot L(j,k) \right), \; i>j
\end{equation}
where $[t]=\{1,...,t\}$. A ``classical'' \chold\ algorithm computes each of these $O(n^3)$ flops, which may be reordered using only commutativity and associativity of addition.  By the no-pivoting and no-distributivity restrictions on $Alg$, when an entry of $L$ is computed, all the entries on which it depends have already been computed and combined by the above formulas, with the sums occurring in any order.
These dependencies form a dependency graph on the entries of $L$, and impose a partial ordering on the computation of the
entries of $L$ (see Figure \ref{fig:dependencies}). That is, when
an entry $L(i,i)$ is computed, by Equation
(\ref{eq:chol-dep-diag}), all the entries $\{L(i,k)\}_{k \in
[i-1]}$ have already been computed.\footnote{While this partial ordering constrains the scheduling of flops, it does not uniquely identify a computation DAG (directed acyclic graph), for the additions within one summation can be in arbitrary order (forming arbitrary subtrees in the computation DAG).} Denote this set of entries by $S_{i,i}$, namely,
\begin{equation}\label{eq:S-diag}
S_{i,i} \equiv \{L(i,k)\}_{k \in [i-1]}
\end{equation}
Similarly, when an entry $L(i,j)$ (for $i>j$) is computed, by
Equation (\ref{eq:chol-dep}), all the entries $\{L(i,k)\}_{k \in
[j-1]}$ and all the entries $\{L(j,k)\}_{k \in [j]}$ have already
been computed. Denote this set by $S_{i,j}$ namely,
\begin{equation}\label{eq:S-lower-tri}
S_{i,j} \equiv \{L(i,k)\}_{k \in  [j-1]} \cup \{L(j,k)\}_{k \in
[j]}
\end{equation}

\begin{lemma}\label{lem:no-masking}
Any ordering of the computation of the elements of $L$ that
respects the partial ordering induced by the above mentioned
directed acyclic graph results in a correct computation of $A \cdot B$.
\end{lemma}

\begin{proof}
We need to confirm that the starred entries $1^*$ and $0^*$ of $T'$
do not somehow ``contaminate'' the desired entries of $L_{32}^T$.
The proof is by induction on the partial order on pairs $(i,j)$
implied by (\ref{eq:S-diag}) and (\ref{eq:S-lower-tri}).
The base case ---the correctness of
computing $L(1,1)$--- is immediate. Assume by induction that all
elements of $S_{i,j}$ are correctly computed and consider the
computation of $L(i,j)$ according to the block in which it resides:

\begin{itemize}
\item
If $L(i,j)$ resides in block $L_{11}$, $L_{21}$ or $L_{31}$ then
$S_{i,j}$ contains only real values, and no arithmetic operations with $0^*$
or $1^*$ occur (recall Figure \ref{fig:dependencies} or Equations
(\ref{eq:GEMM-to-Chol-starred}),(\ref{eq:S-diag}) and
(\ref{eq:S-lower-tri})). Therefore, the correctness follows from
the correctness of the original \chola.

\item
If $L(i,j)$ resides in $L_{22}$ or $L_{33}$ then
$S_{i,j}$ may contain ``starred'' value (elements of $C'$). We
treat separately the case where $L(i,j)$ is on the diagonal and
the case where it is not.

If $i=j$ then by Equation (\ref{eq:chol-dep-diag}) $L(i,i)$ is
determined to be $1^*$ since $T'(i,i)=1^*$ and since adding to,
subtracting from and taking the square root of $1^*$ all result
in $1^*$ (recall Table \ref{tbl:arith} and Equation
(\ref{eq:chol-dep-diag})).

If $i>j$ then by the inductive assumption the divisor $L(j,j)$ of
Equation (\ref{eq:chol-dep}) is correctly computed to be $1^*$
(recall Figure \ref{fig:dependencies} and the definition of $C'$
in Equation (\ref{eq:C})). Therefore, no division by $0^*$ is
performed. Moreover, $T'(i,j)$ is $0^*$. Then $L(i,j)$ is
determined to be the correct value $0^*$, unless $1^*$ is
subtracted (recall Equation (\ref{eq:chol-dep})). However, every
subtracted product (recall Equation (\ref{eq:chol-dep})) is
composed of two factors of the same column but of different rows.
Therefore, by the structure of $C'$, none of them is $1^*$ so
their product is not $1^*$ and the value is computed correctly.

\item
If $L(i,j)$ resides in $L_{32}$ then $S_{i,j}$ may contain
``starred'' values (see Figure \ref{fig:dependencies},
right-hand side, row $j$). However, every subtraction performed
(recall Equation (\ref{eq:chol-dep})) is composed of a product of
two factors, of which one is on the $i$th row (and on a column
$k<j$). Hence, by induction (on $i,j$),  the $(i,k)$ element has
been computed correctly to be  a real value, and by the
multiplication properties so is the product. Therefore no masking
occurs.
\end{itemize}
This completes the proof of Lemma \ref{lem:no-masking}.
\end{proof}

\paragraph{Communication Analysis}
We now know that Algorithm~\ref{alg:Chol-to-MM} correctly multiplies matrices
``classically'', and so has known communication lower bounds
given by
Theorem~\ref{thm:IronyToledoTiskin04}
and
Corollary~\ref{cor:IronyToledoTiskin04}.
It remains to confirm that Step~2 (setting up $T'$)
and Step~4 (returning $L_{32}^T$) do not require much communication,
so that these lower bounds apply to Step 3, running $Alg'$
(recall that Step 1 may be performed off-line and so doesn't count).
Since $Alg'$ is either a small modification of Cholesky to
add ``star'' labels to all data items (at most doubling the bandwidth cost),
or a subset of Cholesky with some operations omitted
(those with starred arguments,
or not leading to the desired output $L_{32}$),
a lower bound on communication for
$Alg'$ is also a lower bound for Cholesky.  \medskip

\begin{theorem}[Main Theorem]
\label{thm:MainTheorem}
Any sequential or
parallel classical algorithm for the \chold\ of $n$-by-$n$
matrices can be transformed into a classical algorithm for
$\frac{n}{3}$-by-$\frac{n}{3}$ matrix-multiplication, in such a
way that the bandwidth cost of the matrix-multiplication algorithm is
at most a constant times the bandwidth cost of the Cholesky algorithm.
\end{theorem}

Therefore any bandwidth or latency cost lower
bound for classical matrix multiplication applies to classical Cholesky,
asymptotically speaking:
\begin{corollary}
\label{cor:lower-sequential}
In the sequential case, with a fast memory of size $M$,
the bandwidth cost lower bound for \chold\ is $\Omega (n^3 / M^{1/2})$,
and
the latency cost lower bound is $\Omega (n^3 / M^{3/2})$.
\end{corollary}

\begin{proof}
Constructing $T'$ (in any data format) requires bandwidth of at
most $18n^2$ (copying a $3n$-by-$3n$ matrix, with another factor
of 2 if each entry has a flag indicating whether it is ``starred''
or not), and extracting $L_{32}^T$ requires another $n^2$ of
bandwidth. Furthermore, we can assume $n^2 < n^3 / M^{1/2}$, i.e.,
that $M < n^2$, i.e., that the matrix is too large to fit entirely
in fast memory (the only case of interest). Thus the bandwidth
lower bound $\Omega (n^3 / M^{1/2} )$ of Algorithm
\ref{alg:Chol-to-MM} dominates the bandwidth costs of Steps 2 and
4, and so must apply to Step 3 (Cholesky). Finally, as each message
delivers at most $M$ words, the latency
lower bound for Step 3 is by a factor of $M$ smaller than its
bandwidth cost lower bound, as desired.
\end{proof}

\begin{corollary}\label{cor:lower:parallel}
In the parallel case, with a 2D layout on $P$ processors, the bandwidth cost lower bound for Cholesky decomposition is
$\Omega (n^2 / P^{1/2})$, and the latency cost lower bound is
$\Omega(P^{1/2})$.
\end{corollary}

\begin{proof}
The argument in the parallel case is analogous to that of
Corollary \ref{cor:lower-sequential}. The construction of input
and retrieval of output at Steps 2 and 4 of Algorithm
\ref{alg:Chol-to-MM} contribute bandwidth of $O\lt(\frac{n^2}{P}
\rt)$. Therefore the lower bound of the bandwidth
$\Omega\lt(\frac{n^3}{P\sqrt{M}}\rt)$ is determined by Step 3, the
Cholesky decomposition. The lower bound on the latency of Step 3
is therefore $\Omega\lt(\frac{n^3}{PM^{3/2}}\rt)$, as each message
delivers at most $M$ words. Plugging in $M=\Theta \lt(
\frac{n^2}{P} \rt)$ yields $B=\Omega(P^{1/2})$.
\end{proof}

\section{Upper Bounds}
\label{sec:upperbound}
In this section we discuss known algorithms for \chold\ and their bandwidth and latency cost analyses, in the sequential two-level memory model, in the sequential hierarchical memory model, and in the parallel computing model.  Throughout this section, we will use the notation $B(n)$ and $L(n)$ for bandwidth and latency costs, respectively, where $n$ is the dimension of the square matrix.  We note that these functions may also depend on the parameters $M$, the size of the fast memory, and $P$, the number of processors, but we omit its reference when it is clear from context.  See also \cite{Bereux08} for latency analysis of two of the \cholas in the two-level (out-of-core) memory model.

In \S 10.1.1 of \cite{Higham96}, standard error analyses of
\chola\ are given. These hold for any ordering of the summation of
Equations (\ref{eq:chol-dep-diag}) and (\ref{eq:chol-dep}), and
therefore apply to all \cholas\ below.

\subsection{Sequential Algorithms}

\subsubsection{Data Storage}
\label{sec:ds}
Before reviewing the various sequential algorithms let us first consider the underlying data-structure which is used for storing the matrix. Although it does not affect the bandwidth cost, it may have a significant influence on the latency cost of the algorithm: it may or may not allow retrieving many relevant words using only a single message.  As a simple example, consider computing the sum of $n \leq M$ numbers in slow memory, which obviously requires reading these $n$ words. If they are in consecutive memory locations, this can be done in one read operation, the minimum possible latency cost. But if they are not consecutive, say they are separated by at least $M-1$ words, this may require $n$ read operations, the maximum possible latency cost.

\begin{figure}[ht]
\begin{center}
\scalebox{.5}{\includegraphics{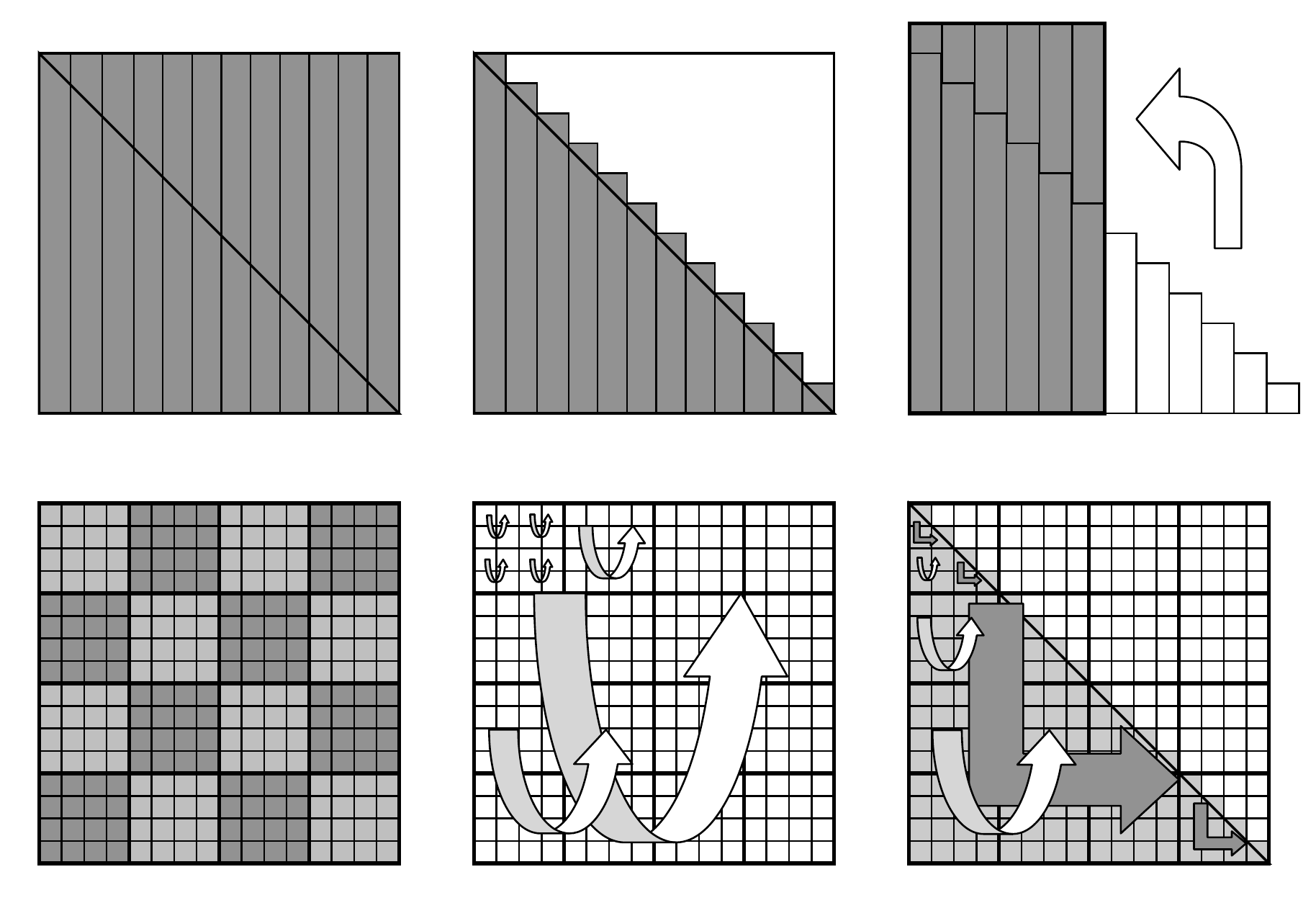}}
\protect\caption{Underlying Data Structures.
\newline Top (column-major): Full, Old Packed, Rectangular Full Packed.
\newline Bottom (block-contiguous): Blocked, Block-Recursive Format, Recursive Full Packed.}
\label{fig:DS}
\end{center}
\end{figure}

The various methods for data storage (only some of which are implemented in \lapack) appear in Figure \ref{fig:DS}. We can partition these data structures into two classes: column-major and block-contiguous.

\paragraph{Column-Major Storage}
The column-major data structures store the matrix entries in column-wise order. This
means that they are most fit for algorithms that access a column
at a time (e.g., the left and right looking na{\"\i}ve algorithms).
However, the latency cost of algorithms that access a block at time
(e.g., \lapack's implementation) will fail to achieve optimality, as retrieving a block of
size $b{\times}b$ requires at least $b$ messages, even if a single
message can deliver $b^2$ words. This means a possible increase in the latency cost by a factor of $b$ (where $b$ is typically the order of $\sqrt{M}$).

As the matrices of interest for \chold\ are symmetric, storing the entire
matrix (known as `Full storage') wastes space. About half of the
space can be saved by using the `old packed' or the
`rectangular full packed' storages. The latter has the advantage
of a more uniform indexing, which allows faster addressing.

\paragraph{Block-Contiguous Storage}
The block-contiguous data structures store the matrix entries in a way that
allows a read or write of a block using a single message. This may
reduce the latency cost of a \chola\ that accesses a $b{\times}b$ block
at a time by a factor of $b$, compared with using column-major storage. We distinguish between two types of block-contiguous storage formats.  The \emph{blocked} data structure stores each
block of the matrix in a contiguous space of the memory. For this,
one has to know in advance the size of blocks to be used by the algorithm (which
is a machine-specific parameter).

The elements of each block may be stored as contiguous sub-blocks, where each sub-block is of a predefined size. The data structure may include several such layers of sub-blocks.  This `layered' data structure may fit the hierarchical memory model (see \S \ref{sec:Hierarchy} for further discussion of this model). The next data structure allows the benefits of block-contiguous storage without knowledge of machine-specific parameters.

The \emph{block-recursive} format \cite{Bader08,ElmrothGustavsonJonssonKagstrom04,FrigoLeisersonProkopRamachandran99,Wise00} (also known as the bit interleaved layout, space-filling curve storage, or Morton ordering format) stores each of the four $n/2{\times}n/2$ submatrices contiguously, and the elements of each submatrix are ordered so that the smaller submatrices are each stored contiguously, and so on recursively. This format is `cache-oblivious'. This means that it
allows access to a single block using a single message (or a
constant number of messages), without knowing in advance the size
of the blocks (see Figure \ref{fig:DS}). Both
the blocked and block-recursive formats have packed versions, where
about half of the space is saved by using the symmetry of the
matrices (see
\cite{ElmrothGustavsonJonssonKagstrom04} for recursive full packed
data structures).  The algorithm in \cite{AndersenGustavsonWasniewski01} uses a hybrid data structure in which only
half the matrix is stored and recursive ordering is used on triangular
sub-matrices and column-major ordering is used on square sub-matrices.

\paragraph{Conversion on the fly}
Since block-contiguous storage yields better latency cost than column-major storage for some algorithms, it can be worthwhile to convert a matrix stored in column-major order to the blocked data structure (with block size $b=\Theta(M)$) before running the algorithm.  This can be done by reading $\Theta(M)$ elements at a time, in a columnwise order (which requires one message) then writing each of these elements to the right location of the new matrix. We write these words using $\Theta(\sqrt{M})$ messages (one per each relevant block). Thus, the total number of messages is $O\lt(\frac{n^2}{\sqrt{M}} \rt) $ which is asymptotically dominated by $O\lt(\frac{n^3}{M^{3/2}}\rt)$ for $M =  \Omega(n)$.  Converting back to column-major storage can be done in similar way (with the same asymptotic cost).

\paragraph{Other Variants of these Data Structures}
Similar to the column-major data structures, there are row-major
data structures that store the matrix entries row-wise. All the
above-mentioned packed data structures have versions that are
indexed to efficiently store the lower (as in Figure \ref{fig:DS})
and upper triangular part of a matrix.

We consider the application of each of the algorithms below to column-major or block-contiguous data structures only.  Adapting each of these algorithms to other compatible data structures (row-major vs. column-major, upper triangular vs. lower triangular, full storage vs. packed storage) is straightforward.

\subsubsection{Arithmetic Count of Sequential Algorithms}
The arithmetic operations of all the sequential algorithms
considered below are exactly the same, up to reordering. The arithmetic
count of all these algorithm is therefore the same, and is given
only once.

By Equations (\ref{eq:chol-dep-diag}),(\ref{eq:chol-dep}) the total arithmetic count is:
$$A(n)=\sum_{i=1}^n \sum_{j=1}^i (2j-1) = \frac{n^3}{3} + \Theta(n^2).$$

\subsubsection{Upper Bound for Matrix Multiplication}
\label{sec:matmul}

In this section we establish the bandwidth and latency costs of matrix multiplication, which is used as a subroutine in some of the \cholas.

We recall the recursive matrix multiplication algorithm of Frigo, Leiserson, Prokop and Ramachandran \cite{FrigoLeisersonProkopRamachandran99}. Their algorithm, given as Algorithm~\ref{alg:RMatMul}, works by a divide-and-conquer approach, where at each step the algorithm splits the largest of three dimensions.

\begin{algorithm}
\protect\caption{$C=RMatMul(A,B)$: A Recursive Matrix Multiplication Algorithm}
\label{alg:RMatMul}
\begin{algorithmic}[1]
\REQUIRE $A,B$, two matrices of dimensions $n{\times}m$ and $m{\times}r$.
\ENSURE $C$, a matrix of dimension $n{\times}r$, so that $C=A \cdot B$
\newline $\setminus \setminus$ {Partition $A,B,C$ by dividing each dimension in half, e.g.
$A=\begin{pmatrix} A_{11} & A_{12} \\ A_{21} & A_{22} \end{pmatrix}$}
\IF{ $n=m=r=1$}
    \STATE $C_{11}=A_{11} \cdot B_{11}$
\ELSIF{$n = \max\{n,m,r \} $}
    \STATE $\begin{pmatrix} C_{11} & C_{12} \end{pmatrix} = RMatMul\lt(\begin{pmatrix} A_{11} & A_{12} \end{pmatrix},B\rt)$
    \STATE $\begin{pmatrix} C_{21} & C_{22} \end{pmatrix} = RMatMul\lt(\begin{pmatrix} A_{21} & A_{22} \end{pmatrix},B\rt)$
\ELSIF {$m=\max \{n,m,r \}$}
    \STATE $C = RMatMul\lt(\begin{pmatrix} A_{11} \\ A_{21} \end{pmatrix}, \begin{pmatrix} B_{11} & B_{12}\end{pmatrix}\rt)$
    \STATE$C = C+RMatMul\lt(\begin{pmatrix} A_{12} \\  A_{22} \end{pmatrix}, \begin{pmatrix} B_{21} & B_{22}\end{pmatrix}\rt)$
\ELSE
    \STATE$\begin{pmatrix} C_{11} \\ C_{21} \end{pmatrix} = RMatMul\lt(A,\begin{pmatrix} B_{11} \\ B_{21}\end{pmatrix}\rt)$
    \STATE$\begin{pmatrix} C_{12} \\ C_{22} \end{pmatrix} = RMatMul \lt(A,\begin{pmatrix} B_{12} \\ B_{22}\end{pmatrix}\rt)$
\ENDIF
\RETURN C
\end{algorithmic}
\end{algorithm}

\begin{lemma}[{\cite{FrigoLeisersonProkopRamachandran99}}]
The bandwidth cost $B_{MM}(n,m,r)$ of multiplying two matrices of dimensions
$n{\times}m$ and $m{\times}r$ is
$$B_{MM}(n,m,r)= \Theta\lt(\frac{nmr}{\sqrt{M}} + nm + mr + nr \rt).$$
\end{lemma}

This result is stated slightly differently in {\cite{FrigoLeisersonProkopRamachandran99}}.  We obtain the form above by setting the cache-line length (or transfer block size) equal to $1$.  The latency cost of this algorithm varies according to the data structure used and the dimensions of the input and output matrices.

\begin{lemma}
When using recursive contiguous-block data structure, the latency cost of matrix multiplication is
$$L_{MM}(n,m,r)= \Theta\lt(\frac{nmr}{M^{3/2}} + \frac{nm+mr+nr}{M} \rt),$$
and when using column-major or row-major layout, the latency cost is
$$L_{MM}(n,m,r)= \Theta\lt(\frac{nmr}{M} + \frac{nm+mr+nr}{\sqrt M} \rt).$$
\end{lemma}
\begin{proof}
Let $m',n',r'$ be the dimensions of the problem the first time all three matrices fit into the fast memory. Then $m',n',r'$ differ by at most a factor of $2$, and thus are all $\Theta(\sqrt{M})$.  Therefore, assuming the recursive contiguous-block data structure, each of the three matrices resides in $O(1)$ square blocks, and so reading and/or writing each matrix incurs a latency cost of $\Theta(1)$. Thus, the total latency cost is $L_{MM}(n)=\Theta\lt(\frac{B_{MM}(n)}{M}\rt)$.

If the column-major or row-major data structures are used, then each such block of size $\Theta(\sqrt{M}){\times}\Theta(\sqrt{M})$ is read or written using $\Theta(\sqrt{M})$ messages (one for each of its rows or columns), and therefore the total latency cost is $L_{MM}(n)=\Theta\lt(\frac{B_{MM}(n)}{\sqrt{M}}\rt)$.

\hfill \end{proof}

We note that Algorithm~\ref{alg:RMatMul} is cache-oblivious.  That is, the algorithm achieves the bandwidth and latency costs given above  with no knowledge of the size of the fast memory.  The iterative blocked matrix multiplication algorithm can also achieve this performance, but the block size must be chosen to be $\Theta\lt(\sqrt M\rt)$.

\subsubsection{Upper Bound for Triangular Solve}
\label{sec:RTRSM}

In this section we establish the bandwidth and latency costs of another subroutine used in some of the \cholas, that of a triangular solve with multiple right hand sides (TRSM).  Algorithm~\ref{alg:rtrsm} is a recursive version of TRSM (a variation of this algorithm appears in \cite{AndersenGustavsonWasniewski01}, designed to utilize a non-recursive, tuned implementation of matrix multiplication).  Below, we assume each matrix multiplication is implemented with the recursive algorithm (Algorithm~\ref{alg:RMatMul}).

\begin{algorithm}
\protect\caption{$X=RTRSM(A,U)$: Recursive Triangular Solver}
\label{alg:rtrsm}
\begin{algorithmic}[1]
\REQUIRE $A,U$ two $n{\times}n$ matrices, $U$ is upper triangular.
\ENSURE $X$, so that $X=A \cdot U^{-1}$
\newline $\setminus \setminus$ {Partition $A,U,X$ by dividing each dimension in half, e.g.
$A=\begin{pmatrix} A_{11} & A_{12} \\ A_{21} & A_{22} \end{pmatrix}$}
\IF{ $n=1$}
    \STATE $X = A / U$
\ELSE
    \STATE $X_{11} = RTRSM(A_{11},U_{11})$
    \STATE $X_{12} = RTRSM(A_{12}-X_{11} \cdot U_{12},U_{22})$
    \STATE $X_{21} = RTRSM(A_{21},U_{11})$
    \STATE $X_{22} = RTRSM(A_{22} - X_{21} \cdot U_{12},U_{22})$
\ENDIF
\RETURN $X$
\end{algorithmic}
\end{algorithm}

\paragraph{Bandwidth Cost}
As no communication is needed for sufficiently small matrices (other than reading the entire input, and writing the output), the bandwidth cost of this algorithm is given by the recurrence
\begin{equation*}
B_{TRSM}(n) = \lt\{ \begin{array}{ll}
4 \cdot B_{TRSM}\lt(\frac{n}{2}\rt)+ 2 \cdot B_{MM}(\frac{n}{2}) & \text{if } n>\sqrt{\frac{M}{3}}  \\
3n^2 & \text{otherwise}
\end{array}\rt.
\end{equation*}
where $B_{MM}(n)=B_{MM}(n,n,n)$ is the bandwidth cost of multiplying two $n$-by-$n$ matrices. Assuming the matrix-multiplication is done using Algorithm~\ref{alg:RMatMul}, we have $B_{MM}(n) = O\lt( \frac{n^3}{\sqrt{M}}+n^2 \rt) $.   The solution to this recurrence (and the total bandwidth cost of recursive TRSM) is then
 $$B_{TRSM}(n)=O\lt(\frac{n^3}{\sqrt M} +n^2\rt).$$

\paragraph{Latency Cost}
Assuming a block-contiguous data structure (recursive, or with the correct block size picked), the latency cost is given by the recurrence
\begin{equation*}
L_{TRSM}(n) \leq \lt\{\begin{array}{ll}
4 \cdot L_{TRSM}\lt(\frac{n}{2}\rt)+ 2\cdot L_{MM}(\frac{n}{2}) & \text{if } n>\sqrt{\frac{M}{3}}  \\
3 & \text{otherwise}
\end{array}\rt.
\end{equation*}
where $L_{MM}(n)=L_{MM}(n,n,n)=O \lt( \frac{n^3}{M^{3/2}} \rt)$ is the latency cost of recursive matrix multiplication.  The solution to this recurrence is
$$L_{TRSM}(n)=O\lt(\frac{n^3}{M^{3/2}}\rt).$$

Again, we note that Algorithm~\ref{alg:rtrsm} is cache-oblivious.  The iterative blocked TRSM algorithm can also achieve this performance, but the block size must be chosen to be $\Theta\lt(\sqrt M\rt)$.

Equipped with the communication costs of matrix multiplication and triangular solve, we proceed to the following subsections, where we present several \cholas\ and provide the communication cost analyses which yield the first five enumerated conclusions in the introduction as well as the results shown in Table~\ref{tbl:sequential}.

\subsubsection{The Na{\"\i}ve Left-Looking Cholesky Algorithm}
\label{sec:naiveleft}
We start with revisiting the two na{\"\i}ve algorithms (left-looking
and right-looking), and see that both have non-optimal bandwidth
and latency, as stated in Conclusion~\ref{con:naive} of the introduction.

\begin{figure}
\begin{center}
\scalebox{.5}{\includegraphics{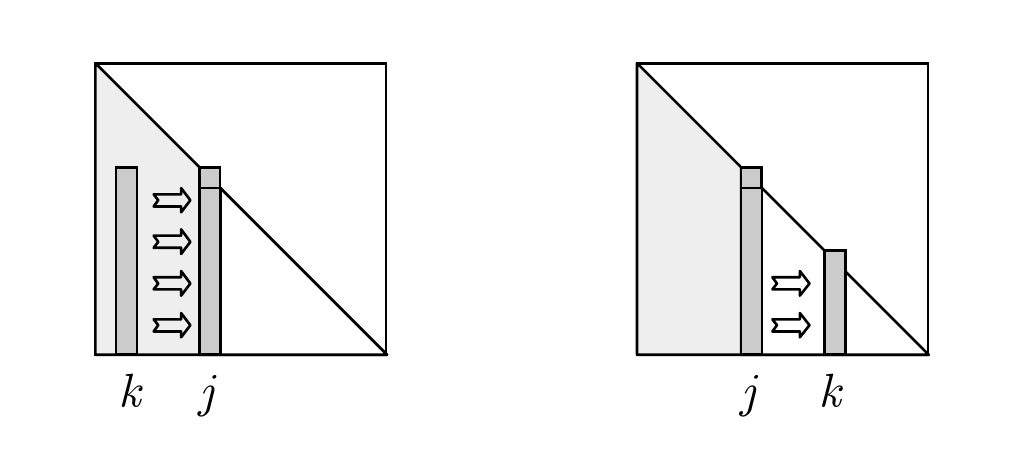}}
\protect\caption{Na{\"\i}ve algorithms.  Left: left-looking. Right: right-looking.  Column $j$ is currently being computed (both algorithms). Light grey is the area already computed.  Column $k$ is the column being read (left-looking) or written (right-looking).}
 \label{fig:Naive}
\end{center}
\end{figure}

The na{\"\i}ve left-looking Cholesky decomposition algorithm is given in Algorithm \ref{alg:naive-left} and Figure \ref{fig:Naive}.  Note that Algorithm \ref{alg:naive-left} assumes that two columns ($k$ and $j$) of the matrix can fit into fast memory simultaneously (i.e., $M>2n$).

\begin{algorithm}
\protect\caption{Na{\"\i}ve left-looking Cholesky
algorithm}
\label{alg:naive-left}
\begin{algorithmic}[1]
\FOR{$j=1$ to $n$}
\STATE read $A(j:n,j)$ from slow memory \label{seqnaiveleft:comm1}
\FOR{$k=1$ to $j-1$}
\STATE read $A(j:n,k)$ from slow memory \label{seqnaiveleft:comm2}
\STATE  update diagonal element: $A(j,j) \leftarrow A(j,j) - A(j,k)^2$
\FOR{$i=j+1$ to $n$}
\STATE update $j^{th}$ column element: $A(i,j) \leftarrow A(i,j)-A(i,k)A(j,k)$
\ENDFOR
\ENDFOR
\STATE calculate final value of diagonal element: $A(j,j) \leftarrow
\sqrt{A(j,j)}$
\FOR{$i=j+1$ to $n$}
\STATE calculate final value of $j^{th}$ column element: $A(i,j) \leftarrow
A(i,j)/A(j,j)$
\ENDFOR
\STATE write
$A(j:n,j)$ to slow memory \label{seqnaiveleft:comm3}
\ENDFOR
\end{algorithmic}
\end{algorithm}

\paragraph{Bandwidth Cost}
Assuming $M>2n$, we follow Algorithm~\ref{alg:naive-left} and communication occurs at lines \ref{seqnaiveleft:comm1},\ref{seqnaiveleft:comm2}, and \ref{seqnaiveleft:comm3}, so the total number of words transferred between fast and slow memory
while executing the entire algorithm is given by
$$\sum_{j=1}^n \lt[ 2(n-j+1) + \lt(\sum_{k=1}^{j-1} (n-j+1)\rt) \rt]=\frac{1}{6}n^3+n^2+\frac{5}{6}n.$$

In the case when $M<2n$, each column $j$ is read into fast memory in segments of size $M/2$.  For each segment of column $j$, the corresponding segments of previous columns $k$ are read into fast memory in sequence to update the current segment.  In this way, the total number of words transferred between fast and slow memory does not change.

\paragraph{Latency Cost}
Assuming $M>2n$ and the matrix is stored in column-major order, each column is contiguous in memory, so the
total number of messages is given by
$$\sum_{j=1}^n \lt[ 2 + \lt(\sum_{k=1}^{j-1} 1\rt) \rt]=\frac{1}{2}n^2+\frac{3}{2}n.$$

In the case when $M<2n$, an entire column cannot be read/written in one message, so the column must be read/written in multiple messages of size $O(M)$.  Again, assuming the matrix is stored in column-major order, segments of columns will be contiguous in memory, and the latency is given by the bandwidth divided by the message size: $O\lt(\frac{n^3}{M}\rt)$.

The algorithm can be adjusted to work one row at a time (up-looking) if the matrix is stored in row-major order (with no change in bandwidth or latency), but a block-contiguous storage format will increase the latency (columns would not be contiguous in memory in this case).  This analysis, along with that of \S~\ref{sec:naiveright}, yields Conclusion~\ref{con:naive} in the introduction.

\subsubsection{The Na{\"\i}ve Right-Looking Cholesky Algorithm}
\label{sec:naiveright}

The na{\"\i}ve right-looking Cholesky decomposition algorithm is given in Algorithm \ref{alg:naive-right} and Figure \ref{fig:Naive}.  Note that Algorithm \ref{alg:naive-right} assumes that two columns  ($k$ and $j$) of the matrix can fit into fast memory simultaneously (i.e., $M>2n$).

\begin{algorithm}
\protect\caption{Na{\"\i}ve right-looking Cholesky algorithm}
\label{alg:naive-right}
\begin{algorithmic}[1]
\FOR{$j=1$ to $n$}
    \STATE read $A(j:n,j)$ from slow memory
    \label{seqnaiveright:comm1}
    \STATE calculate final value for diagonal element: $A(j,j) \leftarrow \sqrt{A(j,j)}$
    \FOR{$i=j+1$ to $n$}
        \STATE calculate final value for $j^{th}$ column element: $A(i,j) \leftarrow A(i,j)/A(j,j)$
    \ENDFOR
    \FOR{$k=j+1$ to $n$}
    \STATE read $A(k:n,k)$ from slow memory
    \label{seqnaiveright:comm2}
        \FOR{$i=k$ to $n$}
            \STATE update $k^{th}$ column element: $A(i,k) \leftarrow A(i,k) - A(i,j)A(k,j)$
        \ENDFOR
        \STATE write $A(k:n,k)$ to slow memory
        \label{seqnaiveright:comm3}
    \ENDFOR
    \STATE write $A(j:n,j)$ to slow memory
    \label{seqnaiveright:comm4}
\ENDFOR
\end{algorithmic}
\end{algorithm}

\paragraph{Bandwidth Cost}
Assuming $M>2n$, we follow Algorithm~\ref{alg:naive-right} and communication occurs at lines \ref{seqnaiveright:comm1}, \ref{seqnaiveright:comm2}, \ref{seqnaiveright:comm3}, and \ref{seqnaiveright:comm4}, so the total number of words transferred between fast and slow memory while executing the entire algorithm is given by
$$\sum_{j=1}^n \lt[ 2(n-j+1) + \sum_{k=1}^{j-1} 2(n-k+1) \rt]= \frac{1}{3}n^3+n^2+\frac{2}{3}n.$$

In the case when $M<2n$, more communication is required.  In order to perform the column factorization, the column is read into fast memory in segments of size $M-1$, updated with the diagonal element (which must remain in fast memory), and written back to slow memory. In order to update the trailing $k^{th}$ column, the $k^{th}$ column is read into fast memory in segments of size $(M-1)/2$ along with the corresponding segment of the $j^{th}$ column.   In order to compute the update for the entire column, the $k^{th}$ element of the $j^{th}$ column must remain in fast memory.  After updating the segment of the $k^{th}$ column, it is written back to slow memory.  Thus, the na{\"\i}ve right-looking algorithm requires more reads from slow memory when two columns of the matrix cannot fit into fast memory, but this increases the bandwidth by only a constant factor.

\paragraph{Latency Cost}
Assuming $M>2n$ and the matrix is stored in column-major order, the total number of messages is given by
$$\sum_{j=1}^n \lt[ 2 + \lt(\sum_{k=j+1}^n 2 \rt)\rt]=n^2+n.$$

As in the case of the na{\"\i}ve left-looking algorithm, when $M<2n$, the size of a message is no longer the entire column.  Assuming the matrix is stored in column-major order, message sizes are of the order equal to the size of fast memory. Thus, for $O(n^3)$ bandwidth, the latency is $O\lt(\frac{n^3}{M}\rt)$.

This right-looking algorithm can be easily transformed into a down-looking row-wise algorithm in order to handle row-major data storages, with no change in bandwidth or latency.  This analysis, along with that of \S~\ref{sec:naiveleft}, yields Conclusion~\ref{con:naive} in the introduction.

\subsubsection{\lapack's {\tt POTRF}}

We next consider an implementation available in \lapack\ (see \cite{LAPACK}) and show that it is bandwidth-optimal when using the right block size (as stated in Conclusion~\ref{con:known-bandwidth} of the introduction) and can also be made latency-optimal, assuming the correct data structure is used (as stated in Conclusion~\ref{con:known-latency} of the introduction).

\begin{figure}
\begin{center}
\scalebox{.5}{\includegraphics{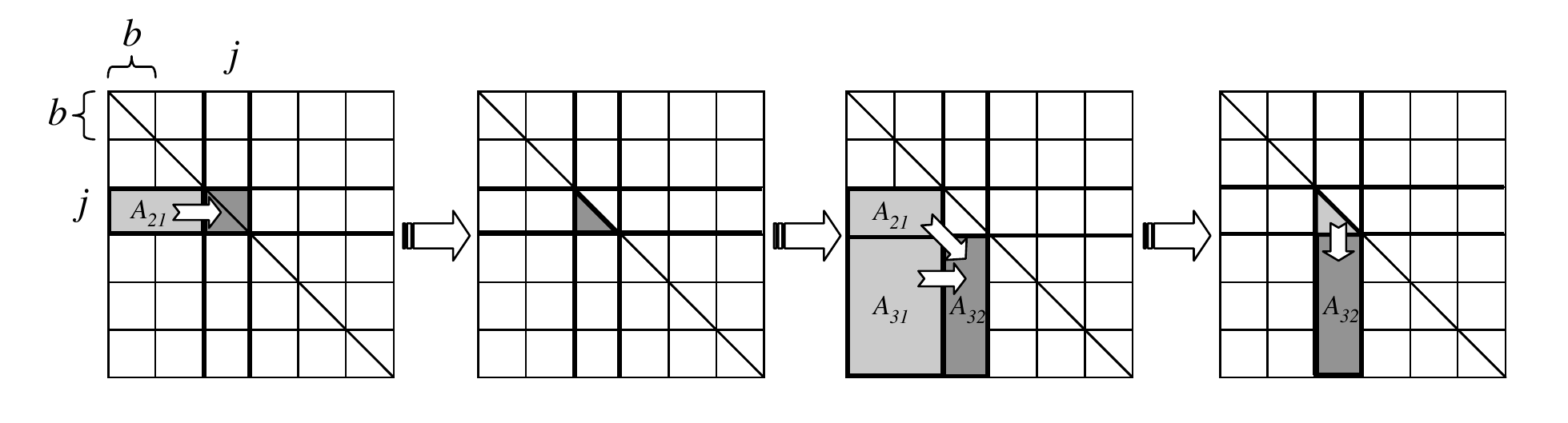}}
\protect\caption{\lapack's {\tt POTRF} Cholesky Decomposition Algorithm}
\label{fig:lapack}
\end{center}
\end{figure}

\lapack's Cholesky algorithm, {\tt POTRF} (Algorithm~\ref{alg:lapack}), is a blocked left-looking algorithm.  For simplicity, we assume the block size $b$ evenly divides the matrix size $n$.  See Figure~\ref{fig:lapack} for a diagram of the algorithm.  Note that by the partitioning of line~\ref{potrf:part}, $A_{21}$ is $b{\times}(j-1)b$, $A_{22}$ is $b{\times}b$, $A_{31}$ is $(n/b-j)b{\times}(j-1)b$, and $A_{32}$ is $(n/b-j)b{\times}b$.

\begin{algorithm}
\protect\caption{LAPACK {\tt POTRF}}
\label{alg:lapack}
\begin{algorithmic}[1]
\label{potrf}
\FOR{$j=1$ to $n/b$}
    \STATE partition matrix so that diagonal block $(j,j)$ is $A_{22}$:$\lt(\begin{matrix}A_{11}&*&*\\A_{21}&A_{22}&*\\A_{31}&A_{32}&A_{33}\end{matrix}\rt)$
    \label{potrf:part}
    \STATE update diagonal block $(j,j)$ ({\tt SYRK}): $A_{22}\leftarrow A_{22}-A_{21}A_{21}^T$
    \label{potrf:syrk}
    \STATE factor diagonal block $(j,j)$ ({\tt POTF2}): $A_{22} \leftarrow \mbox{Chol}(A_{22})$
    \label{potrf:potf2}
    \STATE update column panel ({\tt GEMM}): $A_{32} \leftarrow A_{32}-A_{31}A_{21}^T$
    \label{potrf:gemm}
    \STATE triangular solve for column panel ({\tt TRSM}): $A_{32} \leftarrow A_{32}A_{22}^{-T}$
    \label{potrf:trsm}
\ENDFOR
\end{algorithmic}
\end{algorithm}

The communication costs of Algorithm \ref{potrf} depend on the subroutines which perform symmetric rank-$b$ update, matrix multiply, and triangular solve, respectively. For simplicity, we will assume that the symmetric rank-$b$ update is computed with the same bandwidth and latency as a general matrix multiply.  Although general matrix multiply requires more communication, this assumption will not affect the asymptotic count of the communication required by the algorithm.  We will also assume that the block size $b$ is chosen so that three blocks can fit into fast memory simultaneously; that is, we assume that
$$1\leq b\leq \sqrt{\frac M 3}.$$
In this case, the factorization of the diagonal block (line \ref{potrf:potf2}) requires only $\Theta(b^2)$ words to be transferred since the entire block fits into fast memory. Also, in the triangular solve for the column panel $A_{32}$ (line \ref{potrf:trsm}), the triangular matrix $A_{22}$ is only one block, so the computation can be performed by reading the blocks of $A_{32}$ in sequence and updating them individually in fast memory.  Thus, the amount of communication required by that subroutine is $\Theta(b^2)$ times the number of blocks in the column panel.  Finally, the dimensions of the matrices in the matrix multiply of line \ref{potrf:syrk} during the $j^{th}$ iteration of the loop are $b{\times}b$, $b{\times}(j-1)b$, and $(j-1)b{\times}b$, and the dimensions of the matrices in the matrix multiply of line \ref{potrf:gemm} during the $j^{th}$ iteration are $(n/b-j)b{\times}b$, $(n/b-j)b{\times}(j-1)b$, and $(j-1)b{\times}b$.

\paragraph{Bandwidth Cost}
Under the assumptions above, an upper bound on the number of words transferred between slow and fast memory while executing Algorithm \ref{potrf} is given by
\begin{equation*}
\begin{split}
B(n)\leq\sum_{j=1}^{n/b}\Big[ B_{MM}(b,(j-1)b,b)+\Theta(b^2)+B_{MM}((n/b-j)b,(j-1)b,b) \\
+(n/b-j)\Theta(b^2)\Big]
\end{split}
\end{equation*}
where $B_{MM}(m,n,r)$ is the bandwidth cost required to execute a matrix multiplication of matrices of size $m{\times}n$ and $n{\times}r$ in order to update a matrix of size $n{\times}r$.  Since $B_{MM}$ is nondecreasing in each of its variables, we have
\begin{eqnarray*}
B(n) &\leq& \frac{n}{b}\lt[B_{MM}(b,n,b)+\Theta(b^2)+B_{MM}(n,n,b)+\frac{n}{b}\Theta(b^2)\rt] \\
&\leq& \frac{n}{b}\lt[2B_{MM}(n,n,b)+\Theta(b^2)+\frac{n}{b}\Theta(b^2)\rt].
\end{eqnarray*}
Assuming the matrix--multiply algorithm used by \lapack\ achieves the same bandwidth cost as the one given in \S \ref{sec:matmul} (the iterative blocked algorithm with the correct block size suffices), we have
$$B_{MM}(n,n,b)=\Theta\lt(\frac{n^2b}{\sqrt{M}}+n^2+nb \rt).$$
Since $b\leq \sqrt{M}$ and $b \leq n$, $B_{MM}(n,n,b)=O(n^2)$. Thus,
$$B(n) = O\lt(\frac{n^3}{b} + n^2\rt)$$
and choosing a block size $b=\Theta(\sqrt M)$ gives
$$B(n) = O\lt(\frac{n^3}{\sqrt{M}}+n^2\rt)$$
and achieves the lower bound (as stated in Conclusion~\ref{con:known-bandwidth} of the introduction).  Note that choosing a block size $b=1$ reduces the blocked algorithm to the na{\"\i}ve left-looking algorithm (see Algorithm \ref{alg:naive-left}) which has bandwidth cost $O(n^3)$.

\paragraph{Latency Cost}
Since all reads and writes are done block by block, the optimal latency cost
$$L(n)= O\lt(\frac{n^3}{M^{3/2}} \rt)$$
is achieved if a block-contiguous data structure is used with block size of $\Theta(\sqrt{M})$. However, as the current implementations of \lapack\ use column-major data structures, the latency cost in practice is
$$L(n)= O\lt(\frac{n^3}{M} + \frac{n^2}{\sqrt M} \rt).$$
As mentioned in \S~\ref{sec:ds}, in the case that $M=\Omega(n)$, converting a matrix in column-major order to the blocked layout incurs a bandwidth cost of $O(n^2)$ and latency cost of $O\lt(\frac{n^3}{M^{3/2}}\rt)$.  Thus, in this case, the costs of conversion-on-the-fly and Algorithm~\ref{alg:lapack} combined still achieve the communication lower bounds.  These results are stated in Conclusion~\ref{con:known-latency} of the introduction.

\subsubsection{Rectangular Recursive Cholesky Algorithm}
\label{sec:CholR}

The next recursive algorithm for Cholesky decomposition and its bandwidth analysis follow the recursive $LU$-decomposition algorithm of Toledo (see \cite{Toledo97}). The algorithm given here (Algorithm~\ref{alg:CholR}) is in fact a simplified version of Toledo's: there is no pivoting, and $L=U^T$. For completeness we repeat the bandwidth analysis and provide a latency analysis. The bandwidth proves to be optimal (as stated in Conclusion~\ref{con:known-bandwidth} of the introduction), but the latency does not (as stated in Conclusion~\ref{con:known-latency} of the introduction).

\begin{algorithm}
\protect\caption{$L=RectangularRChol(A)$: Rectangular Recursive Cholesky Algorithm}
\label{alg:CholR}
\begin{algorithmic}[1]
\REQUIRE $A$, an $m{\times}n$ section of positive semidefinite matrix ($m \geq n$). See Figure \ref{fig:CholeskyR} for block partitioning.
\ENSURE $L$, a lower triangular matrix, so that $A=L\cdot L^T$
\IF{ $n=1$}
    \STATE $L=A / \sqrt{A(1,1)}$
\ELSE
    \STATE $\begin{pmatrix} L_{11} \\ L_{21} \\ L_{31} \end{pmatrix} = RectangularRChol\lt(\begin{pmatrix} A_{11} \\ A_{21} \\ A_{31} \end{pmatrix}\rt)$
    \STATE $\begin{pmatrix} A_{22} \\ A_{32} \end{pmatrix} = \begin{pmatrix} A_{22} \\ A_{32} \end{pmatrix} - RMatMul\lt(\begin{pmatrix} L_{21} \\ L_{31} \end{pmatrix}, L_{21}^T \rt)$
    \COMMENT{See \S \ref{sec:matmul} for $RMatMul$.}
    \label{step:recmatmul}
    \STATE $\begin{pmatrix} L_{22} \\ L_{32} \end{pmatrix} = RectangularRChol\lt(\begin{pmatrix} A_{22} \\ A_{32} \end{pmatrix}\rt)$
    \STATE $L_{12} = 0$
\ENDIF
\RETURN L
\end{algorithmic}
\end{algorithm}

\begin{figure}
\begin{center}
\scalebox{.5}{\includegraphics{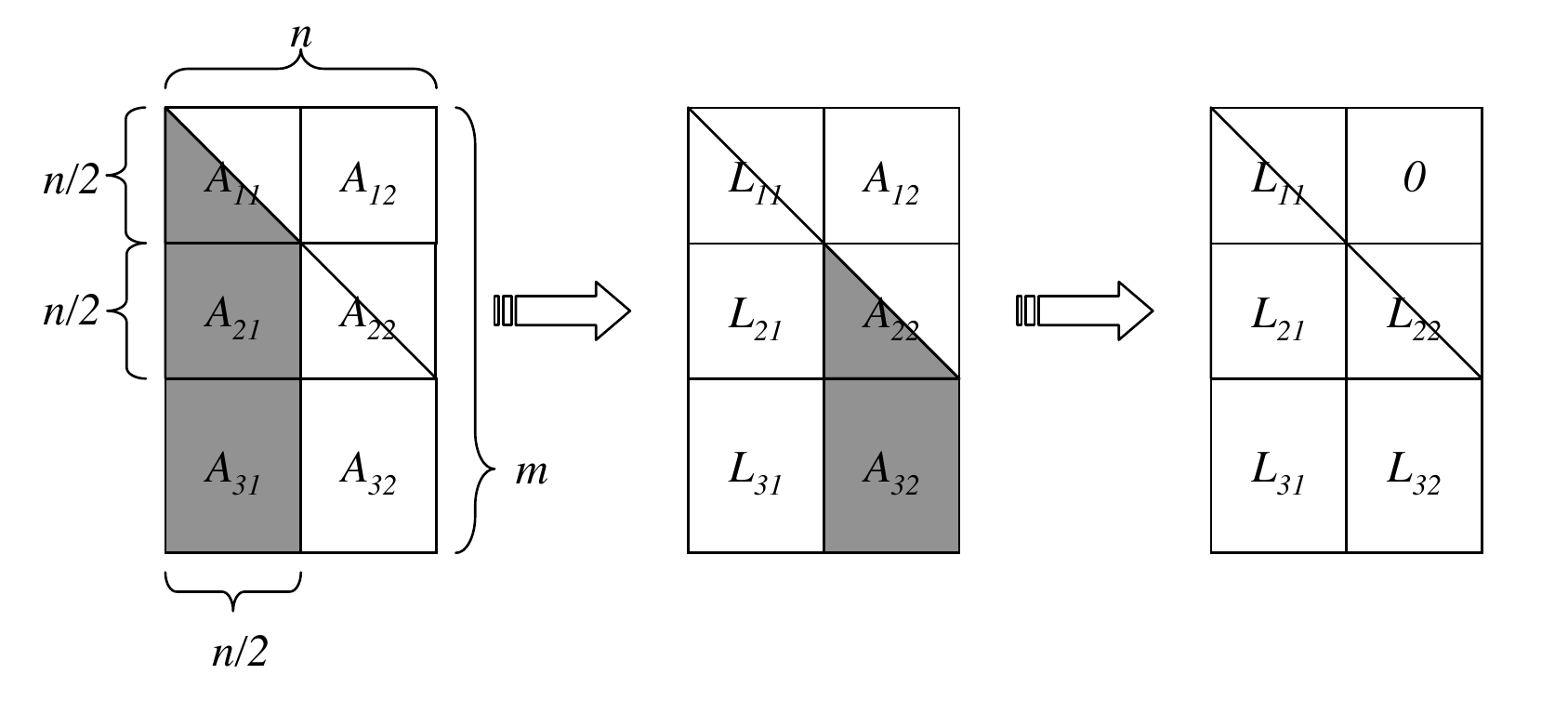}}
\protect\caption{Rectangular Recursive Cholesky Algorithm.}
\label{fig:CholeskyR}
\end{center}
\end{figure}

\paragraph{Bandwidth Cost}
The analysis of the bandwidth cost follows that of Toledo \cite{Toledo97}.  Following Algorithm~\ref{alg:CholR} we either read/write one column (if $n=1$) or make two recursive calls and perform one matrix multiplication. Therefore, the bandwidth cost is given by the recurrence
\begin{equation*}
\label{eqn:TolBW}
B(m,n) = \lt\{\begin{array}{ll}
B\lt(m,\frac{n}{2}\rt)+ B_{MM}(m-\frac n2,\frac n2,\frac n2)+B\lt(m-\frac{n}{2},\frac{n}{2}\rt) & \text{if } n>1  \\
2m & \text{if } n=1
\end{array} \rt. .
\end{equation*}

This recurrence assumes that no more than one column of the matrix can fit in fast memory ($m\leq M$).  We note that Algorithm~\ref{alg:CholR} communicates less than Toledo's algorithm, so the analysis of \cite{Toledo97} provides an upper bound on the solution of this recurrence.  One can use similar analysis to show in the case that $m\geq M$, $B(m,n) = \Omega\lt( \frac{mn^2}{\sqrt{M}}+ mn\log n \rt)$, yielding the following lemma.

\begin{lemma}[{\cite{Toledo97}}]
Given a matrix multiplication subroutine with bandwidth cost equivalent to Algorithm~\ref{alg:RMatMul}, the bandwidth cost of the Rectangular Recursive Cholesky Algorithm on an $m{\times}n$ matrix ($m\geq n$) where at most one column of the matrix can fit into fast memory at a time ($m\geq M$) is
$$B(m,n) = \Theta\lt( \frac{mn^2}{\sqrt{M}}+ mn\log n  \rt).$$
\end{lemma}

Thus, the application of the algorithm to an $n$-by-$n$ matrix (where $n^2>M$) yields a bandwidth cost of
$$B(n,n) = \Theta \left(\frac{n^3}{\sqrt{M}} + n^2\log n \right).$$
This is optimal, provided that $M\leq \frac{n^2}{\log^2 n}$, or, since any algorithm is bandwidth optimal when the whole matrix fits into fast memory, when $n^2 \leq M$.  Hence the range of non-optimality is very small, so we consider Algorithm \ref{alg:CholR} to have optimal bandwidth cost.

We note that in the case where $m < M$, the analysis is slightly different.  In the case that multiple columns fit into fast memory, the bandwidth cost recurrence is
\begin{equation*}
B(m,n) = \lt\{\begin{array}{ll}
B\lt(m,\frac{n}{2}\rt)+ B_{MM}(m-\frac n2,\frac n2,\frac n2)+B\lt(m-\frac{n}{2},\frac{n}{2}\rt) & \text{if } mn>M  \\
\Theta(mn) & \text{if } mn \leq M
\end{array} \rt.
\end{equation*}
where the base case arises when the subproblem fits entirely in fast memory (which may be before Algorithm~\ref{alg:CholR} reaches its base case).  The solution of this recurrence is $B(m,n) = \Theta\lt(\frac{mn^2}{\sqrt M}+mn \log\frac{mn}{M}\rt),$ which is smaller than the cost in the case $m\geq M$.  Because this difference is so small, we omit the distinction in Conclusion~\ref{con:known-bandwidth} and Table~\ref{tbl:sequential}.

\paragraph{Latency Cost}
For the latency cost, we consider the column-major and block-recursive data storage formats.  In order to show that this algorithm does not attain optimal latency cost, we provide lower bounds in each case.

In the case of column-major storage, the latency cost of the entire algorithm is bounded below by the multiplication of Step~\ref{step:recmatmul} at the top of the recursion tree.  This is a square multiplication of size $\frac n2$, so the latency cost with column-major storage is $\Omega\lt(\frac{n^3}{M}\rt)$ (from \S \ref{sec:matmul}). Thus, using column-major storage, Algorithm \ref{alg:CholR} can never achieve optimal latency cost.

In the case of block-recursive storage, the latency of the entire algorithm is bounded below by the latency required to resolve the base cases of the recursion.  Assuming $m>M$, the base case occurs when $n=1$, and the algorithm factors a single column.  However, the block-recursive data structure does not store columns contiguously (up to $2$ elements of the same column could be stored consecutively, depending on the recursive pattern), so reading a column requires $\Omega(n)$ messages.  Since a base case is reached for each column, the latency cost contribution of all the base cases is $\Omega(n^2)$, a lower bound for the total latency cost of Algorithm \ref{alg:CholR}.  Since $n^2$ asymptotically dominates $\frac{n^3}{M^{3/2}}$ for $M>n^{2/3}$, this algorithm cannot achieve optimal latency cost in this case either.

We note that Algorithm~\ref{alg:CholR} is cache-oblivious.

\subsubsection{Square Recursive Cholesky Algorithm}
\label{sec:NewCholR}

The next recursive algorithm is a natural adaptation of the rectangular recursive algorithm of \S~\ref{sec:CholR}.  Instead of dividing the matrix only vertically, because no pivoting is required, we can also divide horizontally, yielding a square recursive algorithm (Algorithm~\ref{alg:NewCholR}).  The algorithm has been considered before, but no asymptotic analysis was given prior to this work.  We believe the algorithm first appears in \cite{Gustavson97} and can also be found in \cite{AhmedPingali00,AndersenGustavsonWasniewski01,GustavsonJonsson01,SimecekTvrdik04}.  Ahmed and Pingali \cite{AhmedPingali00} suggest the algorithm (though without asymptotic analysis of the bandwidth and latency costs). In \cite{SimecekTvrdik04} Simecek and Tvrdik give a detailed probabilistic cache-misses performance; however, no asymptotic claims on the worst-case bandwidth and latency are stated.  We note that the algorithms in \cite{AndersenGustavsonWasniewski01,GustavsonJonsson01} achieve the optimal bandwidth cost, but they assume data storages that prevent achieving optimal latency cost (neither the latency nor the bandwidth cost is explicitly given in those
papers).

The analysis in this section proves that Algorithm~\ref{alg:NewCholR} minimizes the bandwidth cost (as stated in Conclusion~\ref{con:known-bandwidth} in the introduction) and, if the block-recursive data storage is used, the latency cost as well.

\begin{algorithm}
\protect\caption{$L=SquareRChol(A)$: Square Recursive Cholesky
Algorithm}
\label{alg:NewCholR}
\begin{algorithmic}[1]
\REQUIRE $A$, an $n{\times}n$ semidefinite matrix.  See Figure \ref{fig:NewR} for block partitioning.
\ENSURE $L$, a lower triangular matrix, so that $A=L\cdot L^T$
\IF{ $n=1$}
    \STATE $L= \sqrt{A(1,1)}$
\ELSE
    \STATE $L_{11}=SquareRChol(A_{11})$
    \STATE $L_{21}=RTRSM(A_{21}, L_{11}^T)$
        \COMMENT{See \S \ref{sec:RTRSM} for $RTRSM$.}
    \STATE $A_{22}=A_{22}-RMatMul(L_{21}, L_{21}^T)$
        \COMMENT{See \S \ref{sec:matmul} for $RMatMul$.}
    \STATE $L_{22}=SquareRChol(A_{22})$
\ENDIF
\RETURN $L$
\end{algorithmic}
\end{algorithm}

\begin{figure}
\begin{center}
\scalebox{.5}{\includegraphics{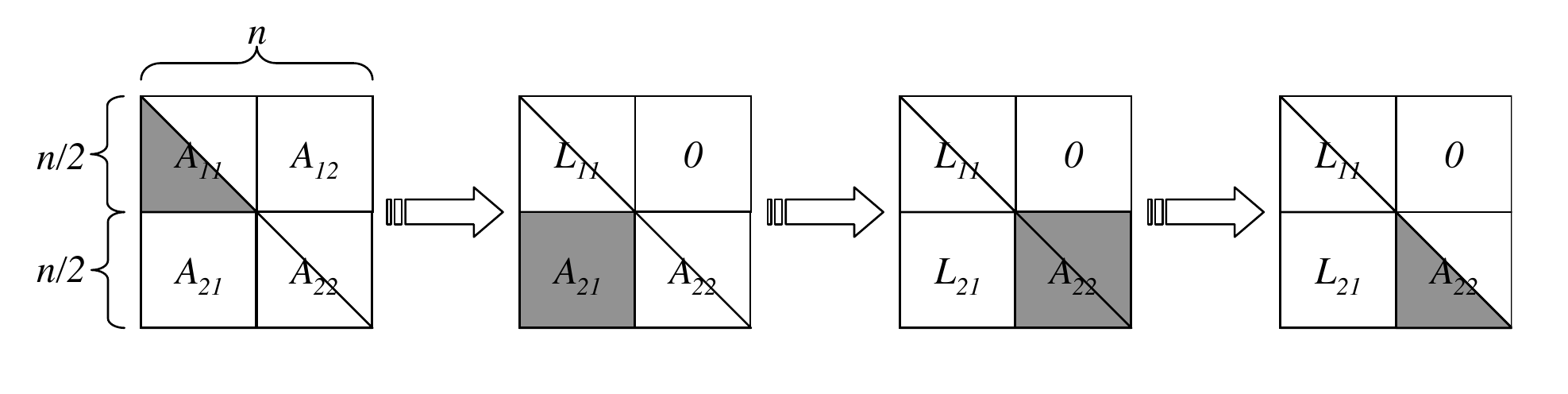}}
\protect\caption{Square Recursive Cholesky Algorithm.}
\label{fig:NewR}
\end{center}
\end{figure}

\paragraph{Bandwidth Cost}
As no communication is needed for sufficiently small matrices (other than reading the entire input and writing the output), the bandwidth cost of this algorithm is given by the recurrence
\begin{equation*}
B(n) = \lt\{\begin{array}{ll}
2 \cdot B\lt(\frac{n}{2}\rt)+ B_{TRSM}\lt(\frac{n}{2}\rt)+B_{MM}\lt(\frac{n}{2}\rt)& \text{if } n>\sqrt{\frac{M}{3}}  \\
2n^2 & \text{otherwise}
\end{array}\rt. .
\end{equation*}
Assuming the matrix multiplication is performed using Algorithm~\ref{alg:RMatMul}, and the triangular solve is performed using Algorithm~\ref{alg:rtrsm}, the solution to this recurrence (and the total bandwidth cost of Algorithm~\ref{alg:NewCholR}) is
 $$B(n)=O\lt(\frac{n^3}{\sqrt M}+n^2\rt).$$

\paragraph{Latency Cost}
The recurrence for the latency cost is similar:
\begin{equation*}
L(n) = \lt\{\begin{array}{ll}
2 \cdot L\lt(\frac{n}{2}\rt)+ L_{TRSM}\lt(\frac{n}{2}\rt)+L_{MM}\lt(\frac{n}{2}\rt)& \text{if } n>\sqrt{\frac{M}{3}}  \\
2 & \text{otherwise, }
\end{array}\rt. .
\end{equation*}
Assuming the matrix is stored using the block-recursive data structure and the recursive subroutines are used, the latency cost is given by
$$L(n)=O\lt(\frac{n^3}{M^{3/2}}\rt).$$

We note that like the Rectangular Recursive Algorithm, Algorithm~\ref{alg:NewCholR} is cache-oblivious.

\subsection{Sequential Algorithm for more than Two Memory Levels}
\label{sec:Hierarchy}
In real computers, there are usually more than two types of memory, and in fact there is a hierarchy of memories, ranging from a very fast small memory to the largest and slowest memory \cite{AndersenGustavsonWasniewski01}. We next consider the lower and upper bounds of \chold\ assuming such a model of {\em hierarchical memory machines}. We observe that none of the known algorithms for \chold\ allows cache-oblivious achievement of optimal latency cost (Conclusion~\ref{con:known-cache-oblivious} of the introduction), except the Square Recursive Cholesky algorithm (Conclusion~\ref{con:new-recursive} of the introduction).

\subsubsection{Lower bound} Explicit lower bounds have been derived for various problems (including matrix multiplication) within the hierarchical memory model (see for example, \cite{Savage95}). Moreover, the known lower bound for
the two-level memory model can be applied here, by considering
any two consecutive levels of the hierarchy as the fast and slow
memories, and treating all faster memories as part of the fast
memory, and all the slower memories as part of the slow one.
Assuming the number of levels is some constant, this may be the
correct lower bound, up to a constant factor. For the \chold, this
gives the following bandwidth and latency cost lower bounds:

\begin{corollary}
Let $Alg$ be a `classical' \chola\ implemented on a machine with some constant
$d$ levels of memory, of size $M_1 \leq \cdots \leq M_d$, with
inverse bandwidth $\beta_1 \leq \cdots \leq \beta_d$ and with
latency $\alpha_1 \leq \cdots \leq \alpha_d$. Then the bandwidth cost
of $Alg$ is
\begin{equation}\label{eq:B-MHM}
\Omega \lt(\sum_{i \in [d-1]} \lt\{ \beta_i \cdot \lt( \frac{n^3}{\sqrt{M_i}}-M_i \rt) \rt\}\rt)
\end{equation}
and the latency cost is
\begin{equation}\label{eq:L-MHM}
\Omega \lt(\sum_{i \in [d-1]} \lt\{ \alpha_i \cdot \frac{n^3}{M_i^{3/2}}\rt\}\rt)
\end{equation}
\end{corollary}

\subsubsection{Upper bounds revisited}
An algorithm may fail to achieve optimal communication costs on a hierarchical memory model even if it achieves optimal performance in the two-level memory model.  For example, an algorithm that includes a parameter (e.g., block size) which allows tuning for better communication performance may be harder or impossible to tune optimally in the hierarchical memory model. On the other hand, as argued in \cite{FrigoLeisersonProkopRamachandran99}, a cache-oblivious algorithm which achieves optimal communication costs in the two-level memory model will also perform optimally on a computer with multiple memory levels.  We next revisit the communication performance of the \cholas\ above in the context of hierarchical memory model.

\paragraph{\lapack's {\tt POTRF}}
In the two-level memory model, the \lapack\ implementation can achieve optimal bandwidth cost, and if the block-contiguous data structure is used, it can also achieve optimal latency cost. To this end one must tune the block size parameter $b$ of the algorithm (and the block size of the data structure).

This has a drawback when applying the \lapack\ algorithm to hierarchical memory machines: setting $b$ to fit a smaller memory level results in inefficient bandwidth and latency costs in the higher levels. Setting $b$ according to a larger memory level results in block I/O that is too large to fit the smaller levels.  Either way, one cannot achieve optimal bandwidth and latency costs between every pair of levels in the memory hierarchy with only one tunable parameter.

\paragraph{Rectangular Recursive Algorithm}
Recall from \S~\ref{sec:CholR} that the recursive \chola\ adopted from \cite{Toledo97} is cache-oblivious and achieves optimal bandwidth cost for practically all ranges of memory sizes, but it cannot guarantee optimal latency cost when $M > n^{2/3}$.  Thus, the algorithm may minimize the bandwidth cost in the hierarchical model, but if, for example, there is a memory level of size greater than $n^{2/3}$ (other than the slowest one holding the original input), then this algorithm yields suboptimal latency cost.

\paragraph{Square Recursive Algorithm}
Recall from \S~\ref{sec:NewCholR} that the Square Recursive algorithm minimizes both bandwidth and latency costs, and it is cache-oblivious.  Thus, it achieves communication optimality in the hierarchical memory model and is the only algorithm to do so.

\subsection{2D Parallel Algorithms}
Let us now consider the so-called $2D$ parallel algorithms, namely those that assume $M=\lt(\frac{n^2}{P}\rt)$ local memory size and start with the $n^2$ matrix elements spread across the processors (i.e, no repetition of elements).

In the parallel model, communication occurs between pairs of processors, and $n$ words sent as one message can be communicated in time $\alpha+\beta n$.  Unlike the sequential model, we do not consider the data structure used to store the parts of the matrix local to each processor.  We assume that the cost of ``packing'' a message into a contiguous block of memory is included in the time $\alpha+\beta n$.  We count the communication cost along the critical path of the algorithm; that is, if many pairs of processors are communicating identical messages simultaneously, the cost along the critical path is that of only one message.

We also consider the cost of a ``broadcast,'' where one processor sends the same message to many other processors.  Broadcasts can be implemented in many different ways, and implementations are often limited by the network topology.  Following \cite{SCALAPACK}, we assume a broadcast of one message incurs a latency cost of $\alpha \cdot \log P$.  That is, if the root processor sends the message to two processors, and each of those processors forward the message to two more, the information can reach $P$ processors after $\log P$ steps.  This implementation requires a topology that includes a binary tree of connections among the $P$ processors (a 2D torus or nearest-neighbor connection topology would not suffice, for example).

Assuming a sufficient network topology, we show upper bounds for bandwidth and latency costs which match the lower bounds up to a logarithmic factor (as stated in Conclusion~\ref{con:parallel} of the introduction).

\subsubsection{The \scalapack\ Implementation}
\label{sec:scalapack}

\begin{figure}
\begin{center}
\scalebox{.5}{\includegraphics{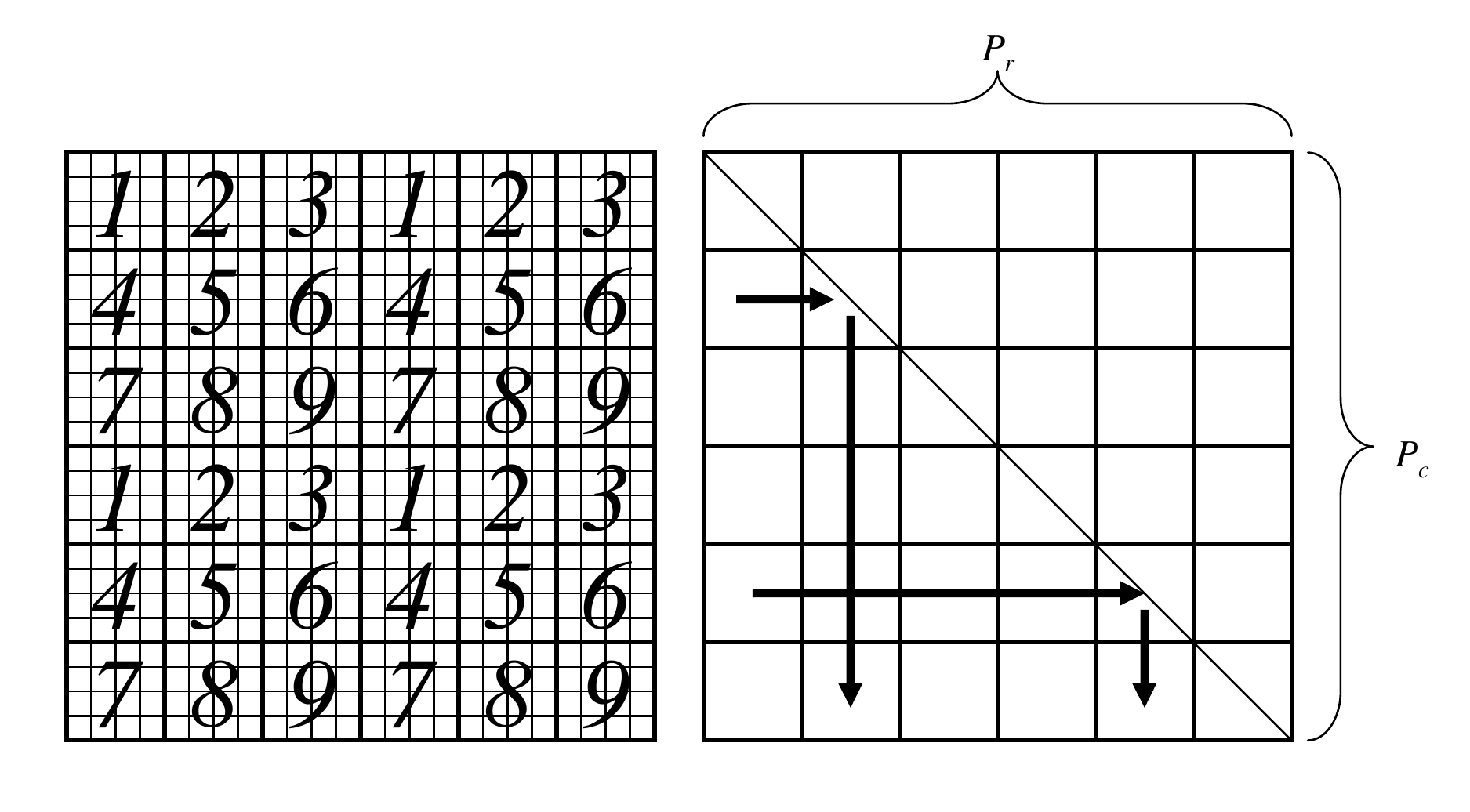}}
\protect\caption{\scalapack's {\tt PxPOTRF}.
\newline
Left: Block-cyclic distribution of the matrix to the processors.  Here $ n=24, b=4, P_c=P_r=3$. \newline
Right: Processor grid: information flow.  Here $P_c=P_r=6$.}
\label{fig:ScaLAPACK}
\end{center}
\end{figure}

The ScaLAPACK \cite{SCALAPACK} routine {\tt PxPOTRF} computes the Cholesky decomposition of a symmetric positive definite matrix $A$ of size $n{\times}n$ distributed over a grid of $P$ processors.  The matrix is distributed block-cyclically, but only half of the matrix is referenced or overwritten (we assume the lower half here).  See Algorithm~\ref{alg:ScaLAPACK} and Figure~\ref{fig:ScaLAPACK}.  We assume a network topology such that each processor column and each processor row are connected via a binary tree.

\begin{algorithm}
\protect\caption{\scalapack\ {\tt PxPOTRF}}
\label{alg:ScaLAPACK}
\begin{algorithmic}[1]
\label{pxpotrf}
\FOR{$j=1$ to $n/b$}
    \STATE processor owning block $(j,j)$ computes Cholesky decomposition
    \STATE broadcast result to processors down processor column
    \FORALL{processor owning blocks in column panel $j$}
        \STATE update blocks in panel with triangular solve
        \STATE broadcast results across processor row
    \ENDFOR
    \FORALL{processor owning diagonal block $(i,i)$ ($i>j$)}
        \STATE re-broadcast results down processor column
    \ENDFOR
    \FORALL{processor owning blocks in trailing matrix}
        \STATE update blocks with symmetric rank-$b$ update
    \ENDFOR
\ENDFOR
\end{algorithmic}
\end{algorithm}

We assume that the block dimension is $b{\times}b$ where $n/b$ is
an integer and the processors are organized in a square grid
($P_r=P_c=\sqrt{P}$) where $n/\sqrt{P}$ is also an integer.  After
computing the Cholesky decomposition of a diagonal block $(j,j)$
locally, the result must be communicated to all processors which
own blocks in the column panel below the diagonal block (i.e.,
blocks $(i,j)$ for $j< i \leq n/b$) in order to update those
blocks.  Since the matrix is distributed block-cyclically, there can
be at most $\sqrt{P}$ such processors.  After the column panel is
updated using a triangular solve,  those results
must be communicated to other processors in order to perform the
rank-$b$ updates to blocks in the trailing matrix. For a given
block $(k,l)$ in the trailing matrix $A_{22}$, the update depends
on the $k^{th}$ block of the column panel (block $(k,j)$) and the
transpose of the $l^{th}$ block of the column panel (block
$(l,j)$).  Thus, after a processor computes an update to block
$(i,j)$ in the column panel, it must broadcast the result to
processors which own blocks in the $i^{th}$ row panel
($P_r=\sqrt{P}$ different processors).  Then, after the processor
owning the diagonal block $(i,i)$ receives that update, it
re-broadcasts to processors which own blocks in the $i^{th}$
column panel ($P_c=\sqrt{P}$ different processors).

The result of the Cholesky decomposition of the diagonal block is
a lower triangular matrix, so the number of words in each message
of the broadcast down the column is only $b(b+1)/2$.  A broadcast
to $P$ processors requires $O\lt(\log P\rt)$ messages.  Any processor
which owns a block in the column panel below the diagonal block
$(j,j)$ will own $\frac{n/b}{\sqrt P}=\frac{n}{b \sqrt P}$ blocks.  Such a
processor computes the updates for all the blocks it owns, and
then broadcasts all the results together (total of
$\frac{nb}{\sqrt P}$ words) to $\sqrt P$ processors (which is done using
$O\lt(\log P\rt)$ messages).  Once a processor owning the corresponding
diagonal blocks receives this message, it re-broadcasts the
message down the column, requiring another $O\lt( \log P \rt)$ messages.
Thus, the total number of messages along the critical path is
$$\sum_{j=1}^{n/b} O\lt(\log P \rt)  = O\lt(\frac{n}{b}\log P\rt)$$
and the number of words along the critical path is
$$\sum_{j=1}^{n/b} \lt[O(b^2\log P) + O\lt(\frac{nb}{\sqrt P}\log P\rt)\rt] = O\lt(\lt(nb+\frac{n^2}{\sqrt{P}}\rt)\log P\rt)$$
Thus, setting the block size at $b=\Theta\lt(\frac{n}{\sqrt{P}}\rt)$, the total
number of messages required is $O\lt(\sqrt{P}\log P\rt)$ and the total number of words is $O\lt(\frac{n^2}{\sqrt{P}} \log P\rt).$  We note that by setting $b=\frac{n}{\sqrt{P}}$, for example, the matrix is no longer stored
block-cyclically.  Instead, each processor owns one large block of the matrix,
and since the full matrix is stored, nearly half of the processors own blocks
which are never referenced or updated.  Further, parallelism is lost as the
algorithm progresses across the column panels.  However, this does not cause an asymptotic degradation in the computation time of the algorithm.  For each of the $\sqrt P$ column panels, there are three phases of computation: Cholesky decomposition of the diagonal block, triangular solve to update the column panel, and matrix multiply to update the trailing matrix.  Thus, setting $b=\frac{n}{\sqrt{P}}$, the computation cost of the algorithm (not including lower order terms) is
$$\sqrt P \lt[ \frac13 \lt(\frac{n}{\sqrt P} \rt)^3 + \lt(\frac{n}{\sqrt P} \rt)^3 + 2\lt(\frac{n}{\sqrt P} \rt)^3\rt] = \frac{10}{3}\cdot\frac{n^3}{P}.$$
Thus, in choosing a large block size to attain the latency cost lower bound, we do not sacrifice the algorithm's ability to meet the computational asymptotic lower bound.\footnote{As noted by Laura Grigori, one can choose a slightly smaller block size and achieve a computational cost of $\frac13\frac{n^3}{P}$ while increasing the latency cost by only a polylogarithmic factor \cite{Grigori09}.}  This analysis yields Conclusion~\ref{con:parallel} in the introduction.

\section{Conclusion}

In this paper we extended known lower bounds on the communication cost (both for bandwidth and for latency)
of conventional ($O(n^3)$) matrix multiplication to Cholesky factorization, and we compared the cost of various Cholesky
decomposition implementations to these lower bounds.  We identified and analyzed existing algorithms that minimized communication for the two-level, hierarchical, and parallel memory models.  We showed that the \scalapack\ implementation minimized communication in the parallel model, up to an $O(\log P)$ factor.  However, the optimal and cache-oblivious sequential algorithm is yet to be implemented in a publicly available library, such as \lapack. Our six main conclusions detailing these results are listed in the introduction.

A `real' computer may be more complicated than any model we have discussed, with both parallelism and multiple levels of memory hierarchy (where each sequential processor making up a parallel computer has multiple levels of cache), or with multiple levels of parallelism (i.e., where each `parallel processor' itself consists of multiple processors), etc. And it may be `heterogenous', with functional units and communication channels of greatly differing speeds. We leave more complicated memory models and lower and upper communication bounds on such processors for future work.

Recently we extended lower bounds to other algorithms in \cite{BDHS10}, including QR decomposition and algorithms for eigenvalue problems, and identified/developed optimal $O(n^3)$ implementations .  We further plan to extend results to Strassen's matrix multiplication and other ``fast'' algorithms.

%\small

\bibliographystyle{siam}
\bibliography{bib_chol_SISC}

\end{document}